\documentclass[12pt]{article}
\usepackage{amssymb,mathrsfs,amsmath}
\usepackage{amscd,amsmath,amsthm,amssymb,amsfonts,enumerate}

\usepackage{graphicx, color}

\usepackage{hyperref}
\hypersetup{
   colorlinks=true,
    linkcolor=blue,
    filecolor=magenta,
   urlcolor=cyan,
    citecolor=blue,
}

\usepackage{subfig}

\usepackage{epsfig} 

\usepackage{amsmath,epsfig,epic,eepic,amsthm,amsfonts}

\newtheorem{thm}{Theorem}[section]

\newtheorem{lem}[thm]{Lemma}
\newtheorem{rem}[thm]{Remark}

\newtheorem{defn}[thm]{Definition}
\newtheorem{exm}[thm]{Example}

\def\ga{\gamma}
\def\P{\mathbb{P}}
\def\E{\mathbb{E}}
\def\Cov{\mathbb{C}{\rm ov}}
\newcommand{\V}{{\mathfrak V}}
\newcommand{\vv}{{\mathfrak v}}

\newcommand{\thmref}[1]{Theorem~{\rm \ref{#1}}}

\newcommand{\disp}{\displaystyle}

\newcommand{\bea}{$$\begin{array}{ll}}
\newcommand{\eea}{\end{array}$$}

\newcommand{\bed}{\begin{displaymath}}
\newcommand{\eed}{\end{displaymath}}
\newcommand{\ad}{&\!\!\!\disp}
\newcommand{\aad}{&\disp}
\newcommand{\barray}{\begin{array}{ll}}
\newcommand{\earray}{\end{array}}
\newcommand{\beq}[1]{\begin{equation} \label{#1}}
\newcommand{\eeq}{\end{equation}}

\newcommand{\bedd}{\bed\begin{array}{l}}
\newcommand{\eedd}{\end{array}\eed}
\newcommand{\al}{\alpha}
\newcommand{\sg}{\sigma}

\newcommand{\e}{\varepsilon}

\newcommand{\nd}{\noindent}

\newcommand{\M}{{\cal{M}}}
\newcommand{\one}{{1}}

\newcommand{\wdt}{\widetilde}

\newcommand{\cd}{(\cdot)}

\def\L{{\cal L}}

\def\({\left(}
\def\]{\right]}
\def\[{\left[}
\def\){\right)}
\def\one{{\hbox{1{\kern -0.35em}1}}}

\makeatletter \@addtoreset{equation}{section}

\def\para#1{\vskip 0.08\baselineskip\noindent{\bf #1}}

\begin{document}

\title{Optimal Control and Numerical Methods for
Hybrid  Stochastic SIS
 Models\thanks{This research was supported in part by Air Force Office of Scientific Research under Grant FA9550-18-1-0268.}
}

\author{Ky Tran\thanks{Department
of Applied Mathematics and Statistics, The
State University of New York in Korea, Yeonsu-Gu, Incheon, Korea 21985, ky.tran@stonybrook.edu.}
\and
George Yin\thanks{Department of Mathematics, University of Connecticut, Storrs, CT
06269, USA. gyin@uconn.edu. }
}


\maketitle

\begin{abstract}	
 This work focuses on  optimal controls of a class of
stochastic SIS epidemic models under regime switching. By assuming that a decision maker can either influence the infectivity period or isolate infected individuals, our aim is to minimize the expected discounted cost due to illness, medical treatment, and the adverse effect on the society. In addition, a model with the incorporation of vaccination is proposed. Numerical schemes are developed by approximating the continuous-time dynamics using Markov chain approximation methods. It is demonstrated that the approximation schemes converge to the
optimal strategy as the mesh size goes to zero.  Numerical examples are provided to illustrate our results.
	
\end{abstract}

\vskip 0.3 true in
\nd{\bf Key words.} Controlled regime-switching diffusion; SIS epidemic model; Markov chain approximation; isolation; vaccination.


\vskip 0.3 true in
\nd{\bf Brief Title.} Control and Numerics for
Hybrid  Stochastic SIS Models

\setlength{\baselineskip}{0.225in}

\newpage

\section{Introduction}\label{sec:int}
  This work focuses on optimal controls and numerical methods for a class of stochastic epidemic models under regime switching.
A SIS system stems from the
compartmental models \cite{KM27}, which subdivides the population into susceptible (S) and infected (I) classes. The word
SIS is
an abbreviation of susceptible-infected-susceptible. That is,
 a susceptible individual becomes infected, and later on becomes susceptible again; see \cite{AM92,Brauer17,KM27} for detailed discussions, and \cite{Gray12} and references therein for different applications.
The models include in particular
diseases being
sexually transmitted and
diseases being transmitted by
bacteria,
in which there is no permanent immunity.
 For such diseases, a promising model is the following classical deterministic SIS epidemic model
	 \beq{ea}
\begin{cases}
dS(t) = \big[\mu N - \beta S(t) I(t) +\gamma I(t) - \mu S(t)\big]dt,\quad S(0)=s_0,\\
dI(t)= \big[\beta S(t)I(t)  - \big( \mu + \gamma\big)I(t)\big]dt, \quad I(0)=i_0,
\end{cases}
\eeq
subject to $S(t)+I(t)=N$, along with the initial values $S(0)=s_0>0$ and $I(0)=i_0>0$, where $S(t)$ and $I(t)$ are  the numbers of susceptible and infected
individuals at time $t$ in a population of size $N$, respectively; $\mu$ and $\gamma^{-1}$ are the average death rate and the average infectious period, respectively. The parameter
$\beta$ is the disease transmission coefficient with $\beta=\lambda/N$ and $\lambda$ being the disease contact rate of an infective individual. More specifically, $\lambda$ is the per day average number of adequate contacts of an infective so that after an adequate contact with an infective, a susceptible individual becomes infected. There is a variety of applications of SIS models. To mention just a few, we refer to \cite{Gray11} for a white noise parameter perturbation of system \eqref{ea}, and \cite{Gray12} for a SIS model with Markovian switching; for SIS models with vaccination we refer to \cite{LM04, LM02, Jiang14, Jiang16}. All the papers above focus on the asymptotic properties of the diseases. The reader can also find related works on long-time behaviors of various epidemic models in \cite{Dang16,NYZ,Ky19}.
	
	 Although SIS epidemic models and
various epidemic models  have been well studied, the work on optimal control of stochastic epidemic models is relatively scarce and the problem is
largely open. The objective of such control problems is to identify effective strategies for minimizing the impacts of infectious diseases through a set of mixed control strategies including treatments, vaccination, isolation, and health promotion campaigns, etc.
Most of the published papers focus on deterministic epidemics and the primary method applied to solve the associated control problems is Pontryagin's maximum principle. To mention just a few, we refer to \cite{Mor74} for a work in the early days of this area, \cite{Beh00, Luca19, Luca17, Clan05, Han11} for a deterministic SIR model with isolation and/or vaccination, \cite{Chen14, Pre13} for controlled epidemic network models. In \cite{Kov18}, the authors studied two stochastic SIS models with external parameter noise and scaled additive noise to minimize the long term average costs. In \cite{Krau18}, the authors investigated SIS meta-population models on networks by comparing different approaches numerically. The recent work \cite{Gran19} studies a SIS model under complete and incomplete observation of the state process by investigating the associated Hamilton-Jacobi-Bellman equation and the Kolmogorov forward equation respectively.
Recently, increasing attention has been devoted to optimal control of the deterministic SIR model due to the appearance of the highly infectious diseases such as COVID-19; see \cite{B20, E20, Le20}.
To the best of our knowledge, not much attention has been given to
stochastic hybrid SIS models
to date. The main difficulties come from the complexity of the model. Even for the stochastic SIS without regime switching, the methods and results in \cite{Gran19} only works under a set of strict assumptions on the model and the cost function.
In this work, we focus on a stochastic SIS model under Markovian switching.
 We develop numerical approximation schemes by approximating the continuous-time dynamics by Markov chains and then showing that the approximations converge to the correct optimal strategy as the mesh size goes to zero.
 Note that the addition of the Markovian switching is to account for environment changes that cannot be modeled by Brownian type of noise, but rather displaying jump behavior.
 We use the Markov chain approximation methodology developed by Kushner and Dupuis \cite{Kushner92}; see also \cite{Song2006}.  Motivated by recent developments in modeling of epidemics, we consider two different models. The first
 one is a direct extension of \eqref{ea}, in which the decision maker is able to control the recovery rate to some extent or isolate a part of the infected individuals. In the second
 model, we incorporate the vaccination into the formulation. That is, the decision maker is able to control the recovery rate and also to take a decision on the fraction of vaccinated individuals.
Moreover, treating a general cost function, we
take into account of any cost associated with the control and either the outbreak size or the infectious burden under the assumption that there are limited control resources.

The rest of the work is organized as follows. Section \ref{sec:for}
begins with the problem formulation.
Section \ref{sec:alg} presents
numerical algorithms based on the Markov chain approximation method.
Section \ref{sec:vac} focuses on a model with the incorporation of vaccination.
In Section \ref{sec:ex}, we present several examples.
Finally, the paper is concluded with some further remarks. To facilitate the reading, all proofs are placed in an appendix at the end of the paper in order not to interrupt the flow of presentation.

\section{Formulation} \label{sec:for}

 Inspired by the recent trend in modeling using regime-switching models, we use a continuous-time Markov chain to model environment changes that are not covered in the usual diffusion models.
Assume throughout the paper that both the Markov chain $\al(t)$
and the scalar standard Brownian motion $w\cd$ are defined
on a complete filtered probability space $(\Omega, \mathcal{F}, \mathcal{F}(t), \P)$, where $\{\mathcal{F}(t)\}$ is a filtration satisfying the usual condition (i.e., right continuous, increasing, and $\mathcal{F}(0)$ containing all the null sets).

Started with
the classical deterministic SIS epidemic model
given in \eqref{ea},
we
normalize the population size to one by replacing $S$ and $I$ with $S/N$ and $I/N$, respectively.
 We obtain
\beq{e.2.2}
\begin{cases}
dS(t) = \big[\mu - \lambda S(t) I(t) +\gamma I(t) - \mu S(t)\big]dt,\quad S(0)=s_0,\\
dI(t)= \big[\lambda S(t)I(t)  - \big( \mu + \gamma\big)I(t)\big]dt, \quad I(0)=i_0,
\end{cases}
\eeq
where $S(t)$ and $I(t)$ are the fraction of susceptible and infected
individuals at time $t$,  respectively.
Throughout the paper, we use $i$ and $s$ to denote the state variables of $I(t)$ and $S(t)$, respectively,  with $i_0$ and $s_0$ being the initial data.
Taking into account the environmental noise, the system parameters $\mu$, $\lambda$, $\gamma$ may experience abrupt changes, which is modeled by a Markov chain $\al\cd$ as in \cite{Gray12}
\beq{e.1}
\begin{cases}
dS(t) = \big[\mu_{\al(t)} - \lambda_{\al(t)} S(t) I(t) +\gamma_{\al(t)} I(t) - \mu_{\al(t)}S(t)\big]dt,\\
dI(t)= \big[\lambda_{\al(t)} S(t)I(t)  - \big( \mu_{\al(t)} + \gamma_{\al(t)}\big)I(t)\big]dt, \\
S(0)=s_0, \quad I(0)= i_0,\quad \al(0)=\ell_0.
\end{cases}
\eeq
We suppose that $\al\cd$ is a continuous time Markov chain taking values in a finite set $\mathcal{M}=\{1, 2, \dots, m_0\}$ generated by $\Lambda=(\Lambda_{\iota \ell})_{m_0\times m_0}$.
 The long time behavior of \eqref{e.1} has been studied in \cite{Gray12}.
A key parameter of the system is the average number of contacts per
infected people
per day
 $\lambda_{\al(t)}$, which can be
 perturbed by white noise; that is, $\lambda_{\al(t)}\to \lambda_{\al(t)} + \sigma_{\al(t)}\dot{w}(t)$ with $w\cd$ being a scalar standard Brownian motion independent of $\al\cd$. Thus,
 \beq{e.1.new}
\begin{cases}
dS(t) = \big[\mu_{\al(t)} - \lambda_{\al(t)} S(t) I(t) +\gamma_{\al(t)} I(t) - \mu_{\al(t)}S(t)\big]dt -\sigma_{\al(t)} S(t)I(t)dw(t),\\
dI(t)= \big[\lambda_{\al(t)} S(t)I(t)  - \big( \mu_{\al(t)} + \gamma_{\al(t)}\big)I(t)\big]dt +\sigma_{\al(t)} S(t)I(t)dw(t), \\
S(0)=s_0, \quad I(0)= i_0,\quad \al(0)=\ell_0.
\end{cases}
\eeq
As in \cite{Gran19, Kov18, Krau18}, we suppose that the natural development of the disease can be influenced by a decision maker. In particular, the decision maker is able to control the magnitude of the recovery rate to some extent. Such a control can be performed by increasing the treatment capacity or the efficiency of medication.
To be more specific, $\gamma_{\al(t)}$ is the recovery rate without any action by the decision maker, while $\gamma_{\al(t)} +C(t)$
is the recovery rate at time $t$ if a control $C(t) \ge 0$ is applied at $t$. As in \cite{Krau18}, $C(t)$ can be interpreted as the product $C(t)=\wdt C(t)
\wdt \zeta$, where $\wdt C(t)$ is the proportion of the infected population treated at time $t$ and $\wdt \zeta$ is the treatment effectiveness.
We can also regard $C(t)$ as per-capital rate of isolation, which
has a direct effect only on infected individual. Thus, the model under consideration can be regarded as a standard isolation model; see \cite{Luca19, Han11}.  To indicate that we have limited resources to handle the epidemic, we suppose that $C\cd$ takes values in a nonempty compact set $\mathcal{U}$ of $[0, \infty)$, where $0\in \mathcal{U}$. The controlled version is given by
\beq{ee2}
\begin{cases}
dS(t) = \big[\mu_{\al(t)} - \lambda_{\al(t)} S(t) I(t) +\gamma_{\al(t)} I(t) - \mu_{\al(t)}S(t)\big]dt -\sigma_{\al(t)} S(t)I(t)dw(t),\\
dI(t)= \big[\lambda_{\al(t)} S(t)I(t)  - \big( \mu_{\al(t)} + \gamma_{\al(t)} + C(t)\big)I(t)\big]dt +\sigma_{\al(t)} S(t)I(t)dw(t), \\
S(0)=s_0, \quad I(0)= i_0,\quad \al(0)=\ell_0.
\end{cases}
\eeq
Before proceeding further, we state the following result regarding the existence of a unique positive solution to \eqref{ee2}.

\begin{thm}\label{thm:one} For any given initial value $(s_0, i_0, \ell_0)\in (0, 1)^2\times \M$ satisfying $s_0+i_0=1$ and $C(t)\equiv c_0\in \mathcal{U}$,  equation \eqref{ee2} has a unique global solution $(S(t), I(t), \al(t))\in (0, 1)^2\times \M$ and $S(t)+I(t)=1$ for any $t\ge 0$
with probability one.
\end{thm}

Given that $I(t)+S(t)=1$,
$I(t)$, the fraction of
 infected individuals, obeys the stochastic Lotka-Volterra model with Markovian switching given by
\beq{ee}dI(t)=I(t)\big[   \lambda_{\al(t)}  -  \mu_{\al(t)} -  \gamma_{\al(t)}  -  \lambda_{\al(t)} I(t) - C(t)\big]dt + \sigma_{\al(t)} I(t)\big[1-I(t)\big]dw(t).\eeq
For convenience, we define
$$b(i, \ell, c) = i\big(   \lambda_{\ell}  -  \mu_{\ell} -  \gamma_{\ell}  -  \lambda_{\ell} i - c\big), \quad a(i, \ell)= \sg^2_\ell i^2(1-i)^2.$$
As in \cite{Luca19, Han11}, we choose the eradication time $\tau$ as the terminal time of planning (that is, time horizon). Formally, let $\zeta$ be a positive constant such that $\xi < 1$, then
\beq{time} \tau=\inf\{t \ge 0:  I(t)\le  \xi\}. \eeq
To make the definition in \eqref{time}  meaningful, we assume that the initial number of infected units $I(0)$ is strictly greater than $\xi$. Thus, $\tau$ is the first time at which the state variable $I$ drops to $\xi$.

Let $\mathcal{A}_{i, \ell
}$ denote the collection of all admissible controls with initial value $(i, \ell)\in (0, 1)\times \mathcal{M}$. A strategy $C\cd$ will be in $\mathcal{A}_{i, \ell
}$ if
$C(t)$ is
$\mathcal{F}(t)$-adapted and $C(t)\in  \mathcal{U}$ for any $t\ge 0$.
We suppose that
increasing the recovery rate is costly and the
treatment of the
infected individuals also create additional costs. The cost is described by the bounded cost function
$F:  [0, 1]\times \M \times [0, \infty)\to [0, \infty)$.
  For a control strategy $C\cd\in \mathcal{A}_{i, \ell}$, we define the cost
  functional as
\beq{e2:10}
J(i, \ell, C\cd):=  \E_{i, \ell}  \int_0^{\tau} e^{-\delta t}  F(I(t), \al(t), C(t)) dt,
\eeq
where $\delta > 0$ is the discounting factor and $\E_{i, \ell}$ denotes the expectation with respect
to the probability law when the process $(I(t), \al(t))$ starts with initial condition $(i, \ell)$. The goal is to minimize the
cost functional and find an optimal  strategy $ C^*\cd$ such that
\beq{e2:11}
J(i,\ell, C^*\cd)=V(i, \ell):= \inf\limits_{C\cd\in \mathcal{A}_{i, \ell}}J(i, \ell, C\cd).
\eeq

Formally, the associated Hamilton-Jacobi-Bellman equation of the underlying problem is given by
 \beq{hjb}
\min\limits_{c\in \mathcal{U}}\Big[b(i, \ell, c)\dfrac{dV(i, \ell)}{d\,i} + \dfrac{1}{2}a(i, \ell)\dfrac{d^2V(i, \ell)}{d\,i^2} + \sum\limits_{\iota=1}^{m_0}\Lambda_{\ell \iota} V(i, \iota) + F(i, \ell, c) -\delta V(i, \ell)\Big]=0,
 \eeq
 for all $(i, \ell)\in [\xi,1]\times \M$, with the boundary condition $V(\xi, \ell)=0$ for any $\ell\in \M$.

\begin{rem} \label{rem:2.1}
{\rm
In the examples in Section \ref{sec:ex}, we will work on a specific cost function of the form
$$F(i, \ell, c) =a_0 + a_1 i + a_2i c^2,$$
where $a_0$, $a_1$, and $a_2$ are the
weighting factors, representing the cost per unit time of the components $1$, $i$, $ic^2$, respectively.
In particular, $\int _0^\tau a_0 dt$ is the cost due to the time period needed for outbreak eradication, while $\int_0^\tau a_1 I(t)dt$ is the cost that infected individual creates for the society due to lost working hours and standard medical care, not including the treatment $C\cd$.  $\int_0^\tau a_2 I(t)C^2(t)dt$ is the cost of treating infected individuals. We assume that the cost is proportional to the number of infected individuals, and the cost per each patient depends quadratically on the treatment effort $C(t)$. As in \cite{Luca19}, if one is interested in minimizing the sum of the eradication time and the total epidemic size, one can use the cost function given by
$$F(i, \ell, c) = a_0 + a_1 \lambda_\ell (1-i)i + a_2ic^2,$$
where $a_0$ and $a_1$ are the weighting factors, representing the cost per unit time of epidemic duration and the cost of a single new infection, respectively.
}
\end{rem}

Our standing assumptions are
as  follows.
\begin{itemize}\item[{\rm (A)}]
\begin{enumerate}
\item
	  The system parameters $\mu_\ell$, $\gamma_\ell$, $\lambda_\ell$ are all nonnegative for each $\ell\in \M$.
	  \item The control set $\mathcal{U}$ is a nonempty compact set in $[0, \infty)$ with $0\in \mathcal{U}$.
	\item  For any $\ell\in \mathcal{M}$, the cost function $F(\cdot, \ell, \cdot)$ is  bounded, continuous, and nonnegative. 	\end{enumerate}
	\end{itemize}

\begin{rem} {\rm
Under assumption (A), for any initial value $(i, \ell)\in  (0, 1)\times \M$ and $C(t)\equiv c\in \mathcal{U}$, equation \eqref{ee} has the unique global solution $(I\cd, \al\cd)\in (0, 1)\times \M$ for all $t\ge 0$; see Theorem \ref{thm:one}.
It follows from the boundedness of  $F\cd$ that $V(i, \ell)<\infty$ for any $(i, \ell)\in (0, 1) \times \M$.

Although our work is motivated by the presence of a Markovian switching random environment, our results and the numerical schemes we developed can also be used to the corresponding control problems of diffusion models without switching. In particular, one can simply take $\M=\{1\}$.
}
 \end{rem}

We are in a position to construct a numerical procedure for
solving the optimal control problem.

\section{Numerical Algorithm}
\label{sec:alg}
Following the Markov chain approximation method in \cite{Kushner92, Song2006}, we  construct a controlled  Markov chain
  in discrete time to approximate the controlled switching diffusions.

\subsection{Approximation Algorithm of the Combined Process}

Let $h>0$ be a discretization parameter for the continuous state variable.
Define
$$S_{h}: = \{i=kh: k\in \mathbb{Z}_+, k\le 1/h\}\times \M.$$
Let $\{(I^h_n, \al^h_n): n\in \mathbb{Z}_+\}$
be a discrete-time controlled Markov chain with state space $S_{h}$
such that
the controlled Markov chain well approximates the local behavior of the controlled diffusion $\big
(I(t), \al(t)\big
)$. At any discrete-time step $n$, the magnitude of the control component $C^h_n$ must be specified. The space of controls is
$\mathcal{U}$.
Let  $C^h=\{C^h_n\}$ be a sequence of controls. We denote by $p^h((i, \ell), (i', \ell')|c)$ the transition probability from state $(i, \ell)$ to another state $(i', \ell')$ under the control $c$. Denote $\mathcal{F}^h_n=\sg\{I^h_k, \al^h_k, C^h_k, k\le n\}$.

The sequence $C^h$ is said to be admissible if it satisfies the following conditions:
\begin{itemize}
	\item[{\rm (a)}]
	$C^h$ is
	$\sigma\{I^h_0, \dots, I^h_{n},\al^h_0, \dots, \al^h_{n},C^h_0, \dots, C^h_{n-1}\}$-{adapted},
	\item[{\rm (b)}]  For any $(i, \ell)\in S_h$, we have
	$$\P\{ (I^h_{n+1}, \al^h_{n+1}) = (i, \ell) | \mathcal{F}^h_n\}= \P\{ (I^h_{n+1}, \al^h_{n+1}) = (i, \ell) | I^h_n, \al^h_n, C^h_n\} = p^h( (I^h_n, \al^h_n), (i, \ell)| C^h_n),$$
	\item[{\rm (c)}] $(I^h_n, \al^h_n)\in S_h$ for all $n\in \mathbb{Z}_+$.
\end{itemize}

The collection
of all admissible control sequences
for initial state $(i, \ell)$ will be denoted by $\mathcal{A}^h_{i, \ell}$.
For each
 $(i, \ell, c)\in S_h\times \mathcal{U}$,
we define
a family of the interpolation intervals $\Delta t^h (i, \ell, c)$. The values of $\Delta t^h (i, \ell, c)$ will be specified later. Then we define
\beq{e.4.3}
\barray
\aad t^h_0 = 0,\quad  \Delta t^h_k = \Delta t^h(I^h_k, \al^h_k, C^h_k),
\quad  t^h_n = \sum\limits_{k=0}^{n-1} \Delta t^h_k.\\
\earray
\eeq
For $(i, \ell)\in S_h$ and $C^h\in \mathcal{A}^h_{i, \ell}$, the 
cost functional for the controlled Markov chain is defined as
\beq{e.4.4}
J^h(i, \ell, C^h) =  \E \sum_{k=1}^{\eta_h} e^{-\delta t_k^h} F(i^h_k, \al^h_k, C^h_k)\Delta t^h_k,
\eeq
with
$$\eta_h=\inf\{n\ge 0: I^h_n\le \xi\}.$$
The value function of the controlled Markov chain is
\beq{e.4.5}
V^h(i, \ell) = \inf\limits_{C^h\in \mathcal{A}^h_{i, \ell}} J^h (i, \ell, C^h).
\eeq
The corresponding dynamic programming equation for the discrete approximation is given by
$$V^h(i, \ell)=\min\limits_{c\in \mathcal{U}} \bigg[e^{-\delta \Delta t^h(i, \ell, c)} \sum\limits_{(i', \ell')\in S_h} V^h(i', \ell')p^h\big( (i, \ell), (i', \ell')| c \big) + F(i, \ell, c)\Delta t^h(i, \ell, c)\bigg],$$
for any $(i, \ell)\in S_h.$

\subsection{Transition Probabilities and Local Consistency}

Let $\E^{h, c}_{i, \ell, n}$, $\Cov^{h, c}_{i, \ell, n}$ denote the conditional expectation and covariance with given
$$\{I_k^h, \al^h_k, C_k^h, k\le n, I_ n^h=i, \al^h_n = \ell, C^h_n=c \},$$
respectively. Define
$$\Delta I_n^h = I_{n+1}^h-I_n^h.$$
Our objective is to define transition probabilities $p^h ((i, \ell), (i', \ell') | c)$ so that the controlled Markov chain $\{(I^h_n, \al^h_n)\}$ is locally consistent with
the controlled switching diffusion \eqref{ee}
in the sense that the following conditions hold:
\beq{e4.1}
\barray
\aad \E^{h, c}_{i, \ell, n}\Delta I_n^h = {a}(i, \ell, c)  \Delta t^h(i, \ell, c) + o(\Delta t^h(i, \ell, c)),\\
\aad \Cov^{h, c}_{i, \ell, n}\Delta I_n^h = a(i, \ell, c)\Delta t^h(i, \ell, c) + o(\Delta t^h(i, \ell, c)),\\
\aad \P^{h, c}_{i, \ell, n}(\al^h_{n+1}=\ell')= \Lambda_{\ell \ell'}\Delta t^h(i, \ell, c)+o(\Delta t^h(i, \ell, c)) \quad \text{for} \quad \ell\ne \ell',\\
\aad \P^{h, c}_{i, \ell, n}(\al^h_{n+1}=\ell)=1+ \Lambda_{\ell\ell}\Delta t^h(i, \ell, c)+o(\Delta t^h(i, \ell, c)),\\
\aad \sup\limits_{n, \ \omega} |\Delta I_n^h| \to 0 \quad \text{as}\quad h \to 0.
\earray
\eeq
Using the procedure in
Inspired by
\cite{Kushner92},
for $(i, \ell)\in S_h$, we define
\beq{e4.2}
\barray
\aad Q_h (i, \ell, c)= a(i, \ell, c) +h|b(i, \ell)| -h^2\Lambda_{\ell\ell}+h,\\
\aad p^h \((i, \ell), (i+h, \ell) |c\) =
\dfrac{a(i, \ell)/2+\big(b(i, \ell, c)\big)^+ h }{Q_h (i, \ell, c)}, \\
\aad p^h \((i, \ell), (i-h, \ell) | c\) =
\dfrac{a(i, \ell)/2+\(b(i, \ell, c)\)^- h}{Q_h (i, \ell, c)}, \\
\aad   p^h \( (i, \ell), (i, \ell') | c\) =\dfrac{h^2 \Lambda_{\ell \ell'} }{ Q_h (i, \ell, c)} \text{ for }\ell\ne \ell' , \quad p^h \( (i, \ell), (i, \ell) | c\) =\dfrac{h  }{ Q_h (i, \ell, c)},\\
\aad \Delta t^h (i, \ell, c)=\dfrac{h^2}{Q_h(i,\ell, c)},
\earray
\eeq
where for a real number $r$,  $r^+=\max\{r, 0\}$,
$r^-=-\min\{0, r\}$.
Set $p^h \((i, \ell), (i', \ell')|c\)=0$ for all unlisted values of $(i', \ell')\in S_h$.
Assumption (A) guarantees that
the transition probabilities in \eqref{e4.2} are well-defined. Using the above transition probabilities,
 we can check that the local-consistency conditions of $\{(I^h_n, \al^h_n)\}$ in \eqref{e4.1} are satisfied.

 \begin{lem}
 The Markov chain $\{(I^h_n, \al^h_n)\}$ with transition probabilities $\{p^h\cd\}$ defined in \eqref{e4.2} is locally consistent with \eqref{ee}.
 \end{lem}

\subsection{Continuous-Time Interpolation and Time Rescaling}

To proceed, we construct a continuous-time interpolation of the approximating chain.
For use in this construction, we define
$n^h(t)=\max\{n: t^h_n\le t\}, t\ge 0$.
The discrete time processes associated with the controlled Markov chain  $\{(I^h_n, \al^h_n)\}$ are defined as follows. Let
\beq{e.4.11}
\barray
\aad
B^{h}_0 = M^{h}_0 =0,\quad  B^{h}_n = \sum\limits_{k=0}^{n-1}
\E^{h}_k \Delta I_k^h,\quad  M^{h}_n  = \sum\limits_{k=0}^{n-1} (\Delta I_k^h -
\E^{h}_k \Delta I_k^h), \quad n\ge 1.
\earray
\eeq
The piecewise constant interpolation processes, denoted by $$I^h\cd,  \al^{h}\cd, B^{h}\cd, M^{h}\cd, C^h\cd$$ are naturally defined as
\beq{e.4.12}
\barray
\aad I^h(t) = I^h_{n^h(t)},\quad \al^h(t) = \al^h_{n^h(t)}, \quad C^h(t) = C^h_{n^h(t)}, \\ \aad  B^{h}(t) = B^{h}_{n^h(t)}, \quad M^{h}(t) = M^{h}_{n^h(t)}, \quad t\ge 0.
\earray
\eeq
Define $\mathcal{F}^h(t)=\sigma\{I^h(s), \al^h(s), C^h(s): s\le t\}$.
We have
\beq{e.4.14} I^h(t)
= i+ B^{h}(t) + M^{h}(t).
\eeq
Recall that $\Delta t^h_k = h^2/Q_h(I^h_k, \al^h_k, C^h_k)$. It follows that
\beq{e.4.15}
\barray
B^{h}(t) \ad = \sum\limits_{k=0}^{n^h(t)-1}
 b (I^h_k, \al^h_k, C^h_k) \Delta t^h_k\\
\ad=\int_0^t  b (I^h(u), \al^h(u), C^h(u)) du-\int_{t^h_{n^h(t)}}^t  b (I^h(u), \al^h(u),C^h(u)) du\\
\ad = \int_0^t  b (I^h(u), \al^h(u), C^h(u) )du + \e^h_1(t),
\earray
\eeq
with $\{\e_1^h\cd\}$
  being an $\mathcal{F}^h(t)$-adapted process satisfying $\lim\limits_{h\to 0} \sup\limits_{t\in [0, T_0]}\E|\e_1^h(t)|=0$  for $T_0\in (0, \infty).$
  Define
  $\tau_h=t^h_{\eta_h}.$
The
cost functional from \eqref{e.4.4} can be rewritten as
\beq{e.4.19}
J^h(i, \ell, C^h)= \E \int_0^{\tau_h} e^{-\delta t} F( I^h(t), \al^h(t), C^h(t)  )dt.
\eeq

\subsection{Relaxed controls}
For our analysis, it is more convenient to use the notion of relaxed controls.
We first briefly recall the notion of ``relaxed control'', which arises naturally in the weak convergence analysis for  the approximation to the optimal control problems.

\begin{defn}{\rm
Let $\mathcal{B}(\mathcal{U}\times [0, \infty))$ be the $\sigma$-algebra of Borel subsets of $\mathcal{U}\times [0, \infty)$. An admissible relaxed control or simply
a relaxed control $m\cd$ is a measure on $\mathcal{B}(\mathcal{U}\times [0, \infty))$ such
that $$m(\mathcal{U}\times [0, t]) = t \quad\text{for all}\quad t\ge 0.$$ Given a relaxed control $m\cd$,
there is a probability measure  $m_t\cd$ defined on the $\sigma$-algebra  $\mathcal{B} (\mathcal{U})$ such that $m(dc dt)=m_t(dc)dt$.

 With the given probability
space, we say that $m\cd$ is an admissible relaxed stochastic control
for $(w\cd, \al\cd)$ or $(m\cd, w\cd, \al\cd)$ is admissible, if (i) for each fixed $t\ge 0$, $m(t, \cdot)$ is a random variable taking values in $\mathcal{R}(\mathcal{U}\times [0, \infty))$, and for each fixed $\omega$, $m(\cdot, \omega)$ is a deterministic relaxed control; (ii) the function defined by $m(A\times [0, t])$ is $\mathcal{F}(t)$-adapted for any $A\in \mathcal{B}(\mathcal{U})$. As a result,
with probability one, there is a measure $m_t(\cdot, \omega)$
on the Borel $\sigma$-algebra $\mathcal{B}(\mathcal{U})$ such that
$m(dcdt) = m_t(dc)dt$. }
\end{defn}

\begin{rem}\label{rem:relax} {\rm
For a sequence  of controls $C^h=\{C^h_n: n\in  \mathbb{Z}_+\}$, we define a sequence of equivalent relaxed controls
as follows.
 First, we set $m_{t^h_n}^h(dc)=\delta_{C^h_n}(dc)$, where $\delta_{C^h_n}\cd$ is the probability measure concentrated at $C^h_n$. Then $m^h\cd$ is defined by $m^h(dcdt)=m_t(dc)dt$. That is,
 $$m^h(B\times [0, t])=\int_0^t \Big(\int_B\delta_{C^h(u)}(dc)\Big)du, \quad B\in \mathcal{B}({\mathcal{U}})\quad \text{and}\quad t\ge 0.$$

Let $\mathcal{R}(\mathcal{U}\times [0, \infty))$ denote the space of all
relaxed controls
on $\mathcal{U}\times [0, \infty)$.
 Then $\mathcal{R}(\mathcal{U}\times [0, \infty))$
 can be metrized using the Prohorov metric in the usual way as in \cite[pp.
 263-264]{Kushner92}. With the Prohorov metric,  $\mathcal{R}(\mathcal{U}\times [0, \infty))$
 is a compact space. It follows that any sequence of relaxed controls has a convergent subsequence.
 Moreover,
a sequence $\eta_n\in \mathcal{R}(\mathcal{U}\times [0, \infty))$ converges to $\eta\in \mathcal{R}(\mathcal{U}\times [0, \infty))$ if and only if for all  continuous functions with compact support $\Psi\cd$ on $\mathcal{U}\times [0, \infty)$,
$$\int_{\mathcal{U}\times [0, \infty)}
\Psi(c, t)\eta_n(dc, dt) \to \int_{\mathcal{U}\times [0, \infty)}
\Psi(c, t)\eta(dc, dt).$$
Note that for a sequence  of ordinary controls $C^h=\{C^h_n:  n\in \mathbb{Z}_+\}$, the associated relaxed control  ${m}^h(dcdt)$ belongs to $\mathcal{R}(\mathcal{U}\times [0, \infty))$.
 Note also that the limits of the ``relaxed
 control representations'' of the ordinary controls might not be ordinary controls, but
 only relaxed controls. }
\end{rem}

With the notion of relaxed control given above, we can write \eqref{e.4.19} as
\beq{e4.12}
J^h(i, \ell, C^h)=J^h(i, \ell, m^h\cd)= \E \int_0^{\tau_h} e^{-\delta t} F( I^h(t),\al^h(t), c) m_t^h(dc)dt.
\eeq
Note also that the value function defined in \eqref{e2:11} can be rewritten as
$$V(i, \ell)=\inf\{ J(i, \ell, m\cd): m\cd \,\,\, \text{ is an admissible relaxed control}\},
$$
where
$$J(i, \ell, m\cd):=  \E_{i, \ell}  \int_0^{\tau} e^{-\delta t} F\big(I(t), \al(t), c\big) m_t(dc)dt.$$

\subsection{Convergence Results}

The proof of the
next lemma can be obtained similar to that of \cite[Theorem 3.1]{Yin03}.

\begin{lem}\label{lem:1}
The process $\{\al^h\cd\}$ converges weakly to $\alpha\cd$, which is a Markov chain with generator $\Lambda=(\Lambda_{\iota \ell})$.
\end{lem}

The following theorems establish the convergence of the approximating Markov chain to the original controlled switching diffusion as well as the convergence of the value functions. To facilitate the reading, all proofs are placed in an appendix at the end of the paper.

\begin{thm}\label{thm:thm}
	Let the chain $\{(I^h_n, \al^h_n) \}$ be constructed with transition probabilities defined in \eqref{e4.2},
	$I^h\cd$, $\al^h\cd$, $B^{h}\cd$, $M^{h}\cd$, $\tau_h\cd$ be the continuous-time interpolation defined in \eqref{e.4.11}-\eqref{e.4.12}, $\{C^h_n\}$ be an admissible control, and $m^h\cd$ be the relaxed control representation of $\{C^h_n\}$.
 Then the following assertions hold.
\begin{itemize}
		\item[\rm (a)]
		The sequence
$${H}^h\cd=\big({I}^h\cd, {\al}^h\cd,   B^{h}\cd,  M^{h}\cd, {m}^h\cd, \tau_h\big)$$
	is tight. 	
	As a result, $({H}^h)_{h>0}$ has a weakly convergent subsequence with limit $${H}\cd=\big({I}\cd, {\al}\cd, {B}\cd, M\cd, {m}\cd, \tau\big).$$
	Moreover, ${I}\cd, {B}\cd, M\cd$  have continuous paths with probability one.
		\item[\rm (b)] For any $t\ge 0$, \beq{e.5.0c}  B(t)=\int_0^t  b ({I}(u),  \al(u), c)m_u(dc)du.
\eeq
		
		\item[\rm (c)]
		${M}\cd$ is a continuous ${\mathcal{F}}(t)$-martingale with quadratic variation $\int_0^t  a( I(u), \al(u) )  du.$  Thus, there is an  $\mathcal{F}(t)$-standard Brownian motion $w(t)$, in which we might have to augment the probability space
such that
\beq{e.5.0f}
M(t)=\int_0^t \sg( I(u), \al(u) )  dw(u).\eeq
			\item[\rm (d)] The limit processes satisfy 	
			\beq{e.5.1}
			 I(t)
= i + B(t) +  M(t), \quad t\ge 0.
			\eeq
	\end{itemize}
\end{thm}

Our main convergence result is given below.

\begin{thm}\label{thm:4.6}  Let $V^h(i, \ell)$ and $V(i, \ell)$ be value functions defined in \eqref{e.4.5} and \eqref{e2:11}, respectively. Then $V^h(i, \ell)\to V(i, \ell)$ as $h\to 0$.
\end{thm}

\section{A Stochastic SIS Model with Vaccination} \label{sec:vac}
Building on our approach of the SIS models,
we consider a SIS epidemic model with vaccination proposed in \cite{LM04, LM02}; see also \cite{Jiang16, Jiang14} for closely related models. We assume that $S(t)$, $I(t)$, and $\V(t)$ are the fraction of infectious, susceptible, and vaccinated individuals at time $t$, respectively. The evolution of $(S\cd, I\cd, \V\cd)$ is given by
\beq{e.4.1}
\begin{cases}
{dS(t)} &= \big[\mu (1-q(t)) -\lambda S(t) I(t) - (\mu + p(t))S(t) + \gamma I(t) + \e \V(t) \big] dt,\\
{d I(t)} &= \big[ \lambda S(t)I(t) - (\mu +\gamma)I(t) \big]dt,\\
{d \V(t)} &=\big[ \mu q(t) +p(t)S(t) - (\mu+\e) \V(t) \big]dt,
\end{cases}
\eeq
where
$q\cd$ is the fraction of vaccinated for newborns,
$p\cd$ is the proportional coefficient of vaccinated for the susceptible, $\e$ is the rate of losing their immunity for vaccinated individuals. The other parameters are understood as in the preceding sections; that is, $\mu$ and $\gamma^{-1}$ are the average death rate and the average infectious period respectively,  $\lambda$ is the disease contact rate of an infective individual. Note that if $\V(0)=p(t)=q(t)=0$ for $t\ge 0$, then $\V(t)=0$ for any $t\ge 0$ and the system \eqref{e.4.1} reduces to \eqref{e.2.2}.

 Using the same methods as in the preceding sections, we obtain the following controlled SIS system under regime switching in random environment
\beq{ee.3}
\begin{cases}
{dS(t)} &= \big[\mu_{\al(t)} (1-q(t)) -\lambda_{\al(t)} S(t) I(t) - (\mu_{\al(t)} + p(t))S(t) + \gamma_{\al(t)} I(t) + \e_{\al(t)} \V(t) \big] dt \\
& \hspace{5.5cm} - \sigma_{\al(t)} S(t)I(t)dw(t),\\
{d I(t)} &= \big[ \lambda_{\al(t)} S(t)I(t) - \big(\mu_{\al(t)} +\gamma_{\al(t)} +C(t)\big)I(t) \big]dt + \sigma_{\al(t)} S(t)I(t)dw(t),\\
{d \V(t)} &=\big[ \mu_{\al(t)}q(t) +p(t)S(t) - (\mu_{\al(t)}+\e_{\al(t)}) \V(t) \big]dt,
\end{cases}
\eeq
where $\al\cd$ is a Markov chain taking values in $\mathcal{M}=\{1, 2, \dots, m_0\}$ and $w\cd$ is a scalar standard Brownian motion independent of $\al\cd$.
We suppose that $C\cd$ takes values in a nonempty compact set $\mathcal{U}$ of $[0, \infty)$.  Suppose $\mathcal{V}_p$ and $\mathcal{V}_q$ are nonempty compact subsets of $[0, 1]$. Moreover, $p\cd$ and $q\cd$ takes the values in $\mathcal{V}_p$ and $\mathcal{V}_q$, respectively.
Compared to the related formulations in \cite{Luca19, Han11}, we also consider the possibility that vaccinated individuals lose their immunity.
We suppose that the initial value is $(s_0, i_0, \vv_0)\in (0, 1)^2\times [0, 1)$ satisfying $s_0+ i_0+\vv_0=1$. Thus, we exclude the trivial case $i_0=0$ or unrealistic cases such as $i_0=1$ and $s_0=0$.  A fundamental result on
the existence of a unique global solution of \eqref{ee.3} is given below.

\begin{thm}\label{thm:two} For any given initial value $(s_0, i_0, \vv_0, \ell_0)\in (0, 1)^2\times [0,1)\times \M$ satisfying $s_0+i_0+\vv_0=1$, $C(t)\equiv c_0\in \mathcal{U}$, $p(t)\equiv p_0\in \mathcal{V}_p$ and $q(t)\equiv q_0\in \mathcal{V}_q$,
 the equation \eqref{ee.3} has a unique global solution $(S(t), I(t), \V(t), \al(t))\in (0, 1)^2\times [0, 1)\times \M$ and $S(t) + I(t)+\V(t)=1$
  for any $t\ge 0$
with probability one.
\end{thm}

Since $S(t) +I(t)+\V(t)=1$ for any $t\ge 0$, we need only
consider the last two equations; that is,
\beq{exh}
\begin{cases}
{d I(t)} &= \big[ \lambda_{\al(t)} \big(1-I(t)-\V(t)\big)I(t) - \big(\mu_{\al(t)} +\gamma_{\al(t)} +C(t)\big)I(t) \big]dt\\
&\qquad\qquad + \sigma_{\al(t)} \big(1-I(t)-\V(t)\big)I(t)dw(t),\\
{d \V(t)} &=\big[ \mu_{\al(t)}q(t) +p(t)\big(1-I(t)-V(t)\big) - (\mu_{\al(t)}+\e_{\al(t)}) \V(t) \big]dt.
\end{cases}
\eeq

Let $\mathcal{A}_{i, \vv, \ell
}$ denote the collection of all admissible controls with initial value $(i, \vv, \ell)\in (0, 1)\times [0, 1)\times \mathcal{M}$ satisfying $0<i+\vv\le 1$; that is,
$$I(0)=i, \quad \V(0)=\vv.$$
 Then Theorem \ref{thm:two} indicates that $\mathcal{A}_{i, \vv, \ell
}\ne \emptyset.$
The control component is $(C\cd, p\cd, q\cd)$.
A strategy $(C\cd, p\cd, q\cd)$ will be in $\mathcal{A}_{i, v, \ell
}$ if $\big(C(t), p(t), q(t)\big)$ is $\mathcal{F}(t)$-adapted and $$\big(C(t), p(t), q(t)\big)\in  \mathcal{U}\times \mathcal{V}_p\times \mathcal{V}_q \ \hbox{ for any  } \ t\ge 0.$$
The cost is described by the bounded cost function
$F:  [0, 1]\times \M \times \mathcal{U}\times \mathcal{V}_p\times \mathcal{V}_q\to [0, \infty)$.
  For a control strategy $(C\cd, p\cd, q\cd)\in \mathcal{A}_{i, \vv, \al}$, we define the
cost  functional as
\beq{e2:10}
J(i, \vv, \ell, C\cd, p\cd, q\cd):=  \E_{i, \vv, \ell}  \int_0^{\tau} e^{-\delta t}  F\big(I(t), \al(t), C(t), p(t), q(t)\big) dt,
\eeq
where $\delta > 0$ is the discounting factor and $\E_{i, v,\ell}$ denotes the expectation with respect
to the probability law when the process $(I(t), \V(t), \al(t))$ starts with initial condition $(i, \vv, \ell)$. The goal is to minimize the 
cost functional and find an optimal  strategy $(C^*\cd, p^*\cd, q^*\cd)$ such that
\beq{e2:11.3}
J(i,\vv, \ell, C^*\cd, p^*\cd, q^*\cd)=V(i, \vv, \ell):= \inf\limits_{(C\cd, p\cd, q\cd)\in \mathcal{A}_{i, \vv, \al}}J(i, \vv, \ell, C\cd, p\cd, q\cd).
\eeq
We use the same method as in the preceding section to construct a discrete-time controlled Markov chain $\{(I^h_n, \V^h_n, \al^h_n): n\in \mathbb{Z}_+\}$ with the state space
$$T_{h}: = \big(\{(k_1h, k_2h): k\in \mathbb{Z}_+, k_1+k_2\le 1/h\}\cap [0, 1]\big)^2\times \M.$$
At any discrete-time step $n$, the magnitude of the control component $(C^h_n, p^h_n, q^h_n)$ must be specified. The space of controls is
$\mathcal{U}\times \mathcal{V}_p\times \mathcal{V}_q$.
Let  $(C^h, p^h, q^h)=\{(C^h_n, p^h_n, q^h_n)\}$ be a sequence of controls.
Denote $\mathcal{F}^h_n=\sg\{I^h_k, \V^h_k,  \al^h_k, C^h_k, p^h_k, q^h_k, k\le n\}$.

The sequence $(C^h, p^h, q^h)$ is said to be admissible if it satisfies the following conditions:
\begin{itemize}
	\item[{\rm (a)}]
	$(C^h, p^h, q^h)$ is
	$\sigma\{I^h_0, \V^h_0, \al^h_0\dots, I^h_{n}, \V^h_n, \al^h_{n}, C^h_0, p^h_0, q^h_n\dots, C^h_{n-1}, p^h_{n-1}, q^h_{n-1}\}$-{adapted},
	\item[{\rm (b)}]  For any $(i, \vv, \ell)\in T_h$, we have
	\begin{equation*}
	\barray
	&\P\{ (I^h_{n+1}, \V^h_{n+1}, \al^h_{n+1}) = (i, \vv, \ell) | \mathcal{F}^h_n\}\\
	&\qquad = \P\{ (I^h_{n+1}, \V^h_{n+1}, \al^h_{n+1}) = (i, \vv, \ell) | I^h_n, \V^h_n, \al^h_n, C^h_n, p^h_n, q^h_n\} \\
	&\qquad = p^h( (I^h_n, \V^h_n, \al^h_n), (i,\vv, \ell)| C^h_n, p^h_n, q^h_n),
	\earray
	\end{equation*}
	\item[{\rm (c)}] $(I^h_n, \V^h_n, \al^h_n)\in T_h$ for all $n\in \mathbb{Z}_+$.
	\end{itemize}

The class of all admissible control sequences $(C^h, p^h, q^h)$ for initial state $(i, \vv, \ell)$ will be denoted by $\mathcal{A}^h_{i, \vv, \ell}$.
We need to define the transition probabilities $p^h((i, \vv, \ell), (i', \vv', \ell')|c, p, q)$ so that the controlled Markov chain $\{(I^h_n, \V^h_n,  \al^h_n): n\in \mathbb{Z}_+\}$ is locally consistent with respect to the controlled diffusion \eqref{ee.3}. To proceed, we denote
\beq{}
\barray
\aad b_1(i, \vv, \ell, c) = \lambda_\ell (1-i-\vv)i - (\mu_\ell +\gamma_\ell + c)i, \quad a_1(i, \vv, \ell) = \sg_\ell^2 (1-i-\vv)^2i^2,\\
\aad b_2(i, \vv, \ell, c, p, q) = \mu_\ell q + p(1-i-\vv) - (\mu_\ell +\e_\ell)\vv.
\earray
\eeq
In particular, for $(i, \vv, \ell)\in T_h$ and $(c, p, q)\in \mathcal{U}\times \mathcal{V}_p\times \mathcal{V}_q$, we define
\beq{e4.2}
\barray
\aad Q_h (i, \vv, \ell, c, p, q)= a_1(i, \vv, \ell) +h|b_1(i, \vv, \ell, c)| + h|b_2(i, \vv, \ell, c, p, q) | -h^2\Lambda_{\ell \ell}+h,\\
\aad p^h \((i, \vv, \ell), (i+h, \vv, \ell) |c, p, q\) =
\dfrac{a_1(i, \vv, \ell)/2+\big(b_1(i, \vv, \ell, c)\big)^+ h }{Q_h (i, \vv, \ell, c, p, q)}, \\
\aad p^h \((i, \vv, \ell), (i-h, \vv, \ell) | c, p, q\) =
\dfrac{a_1(i, \vv, \ell)/2+\(b_1(i, \vv, \ell, c)\)^- h}{Q_h (i, \vv, \ell, c, p, q)}, \\
\aad p^h \((i, \vv, \ell), (i, \vv+h, \ell) |c, p, q\) =
\dfrac{\big(b_2(i, \vv, \ell, c, p, q)\big)^+ h }{Q_h (i, \vv, \ell, c, p, q)}, \\
\aad p^h \((i, \vv, \ell), (i, \vv-h, \ell) | c, p, q\) =
\dfrac{\(b_2(i, \vv, \ell, c, p, q)\)^- h}{Q_h (i, \vv, \ell, c, p, q)}, \\
\aad   p^h \( (i, \vv, \ell), (i, \vv, \ell') | c, p, q\) =\dfrac{h^2 \Lambda_{\ell \ell'} }{ Q_h (i, \vv, \ell, c, p, q)} \text{ for } \ell \ne \ell',\\
\aad p^h \( (i, \vv, \ell), (i, \vv, \ell) | c, p, q\) =\dfrac{h  }{ Q_h (i, \vv, \ell, c, p, q)},\quad \Delta t^h (i, \vv, \ell, c, p, q)=\dfrac{h^2}{Q_h(i, \vv, \ell, c, p, q)}.
\earray
\eeq
For $(i, \vv, \ell)\in T_h$ and $(C^h, p^h, q^h)\in \mathcal{A}^h_{i, \vv, \ell}$, the cost functional
 for the controlled Markov chain is given by
\beq{e.4.4.4}
J^h(i, \vv, \ell, C^h, p^h, q^h) =  \E \sum_{k=1}^{\eta_h} e^{-\delta t_k^h} F(i^h_k, \al^h_k, C^h_k, p^h_k, q^h_k)\Delta t^h_k.
\eeq
The value function is
\beq{e.4.4.5}
V^h(i, \vv, \ell) = \inf\limits_{(C^h, p^h, q^h)\in \mathcal{A}^h_{i, \vv, \ell}} J^h (i, \vv,\ell, C^h, p^h, q^h).
\eeq
The main convergence result in this case is given below.

\begin{thm}\label{thm:4.6.3} Let $V(i, \vv, \ell)$ and $V^h(i, \vv, \ell)$ be the value functions defined in \eqref{e2:11.3} and \eqref{e.4.4.5}, respectively. Then $V^h(i, \vv, \ell)\to V(i, \vv, \ell)$ as $h\to 0$.
\end{thm}

\section{Examples} \label{sec:ex}

Throughout this section, we suppose the discounting factor is $\delta = 0.05$. Also, we  work with $\M= \{1, 2\}$ and $\xi=0.02$.

\begin{exm} \label{ex1}
{\rm
We start with the model given by \eqref{ee}. That is,
\beq{eee}dI(t)=I(t)\big[   \lambda_{\al(t)}  -  \mu_{\al(t)} -  \gamma_{\al(t)}  -  \lambda_{\al(t)} I(t) - C(t)\big]dt + \sigma_{\al(t)} I(t)\big[1-I(t)\big]dw(t).\eeq
For $(i, \ell)\in (\xi, 1]\times \M$, we take the initial control $C_0(i, \ell)\equiv \max \mathcal{U}$
and set the initial value
 $V^h_0(i, \ell)\equiv \int_0^\infty e^{-\delta t}F(1, \ell, C_0) dt$.
We  outline how to find the sequence $V^h_n(\cdot)$ as follows.
For each $(i, \ell)\in S_h$ and control $c\in \mathcal{U}$,
 we compute
$$V_{n+1}^h(i, \ell \,|\, c) =e^{-\delta \Delta t^h(i, \ell, c) }\sum\limits_{(i', \ell')\in S_h} V^{h}_n (i', \ell'
) p^h \big((i, \ell), (i', \ell') | c \big) + F(i, \ell, c)\Delta t^h(i, \ell, c).$$ Note also that if $i\le \xi$, then  $V_n^h(i, \ell)=0$ for any $\ell$ and $n$.
We choose the control $C^h_{n+1}(i, \ell)$ and record an improved value
$V^h_{n+1}(i, \ell)$ by
$$C^h_{n+1}(i, \ell) = {\rm arcmin}_{c\in \mathcal{U}} V_{n+1}^h(i, \ell \,|\, c),  \quad V^h_{n+1} (i, \ell) = V^h_{n+1}(i, \ell\,|\, C^h_{n+1}(i, \ell)).$$
The iterations stop as soon as the
increment
$V^h_{n}\cd-V^h_{n+1}\cd$
reaches a predetermined
tolerance level.
We set the error tolerance to be $10^{-8}$.
Our specific example is motivated by \cite[Example 6.2.1]{Gray12}. The system parameters are given by
$$\mu_1=0.45, \quad \mu_2=0.05, \quad \gamma_1=0.35, \quad \gamma_2 = 0.15, \quad \lambda_1 = 2, \quad \lambda_2=2.4, \quad \sg_1=0, \quad \sg_2=1.$$
Suppose that the generator $\Gamma$ of the Markov chain $\al\cd$ is given by $$\Lambda_{11}=-1, \quad \Lambda_{12}=1, \quad \Lambda_{21}=1, \quad \Lambda_{22}=-1.$$
  The  set of controls is given by
$\mathcal{U}=\{k/5: 0\le k\le 15\}.$ Thus, there are 16 control levels.

In the first example,
we consider the cost function
\beq{2e:ee} F(i, \ell, c)=1+\ell i+\ell ic^2;\eeq
see Remark 2.1 for its interpretation. The value function and optimal control is shown in Figure \ref{fig1A}. In each regime $\ell$,  there is a level $L(\ell)$ such that the optimal control is the maximum control effort for $i<L(\ell)$.
It appears that the value function in regime 2 is higher than that in regime 1, while the optimal control in regime 2 is strictly smaller than that in regime 1 possibly because of the higher cost function in regime 2; that is $F(i, 2, c)>F(i, 1, c)$ for any $(i, c) \in (0, 1)\times \mathcal{U}$.

\begin{figure}[h!tb]
	\begin{center}
		\includegraphics[height=3.5in,width=5.5in]{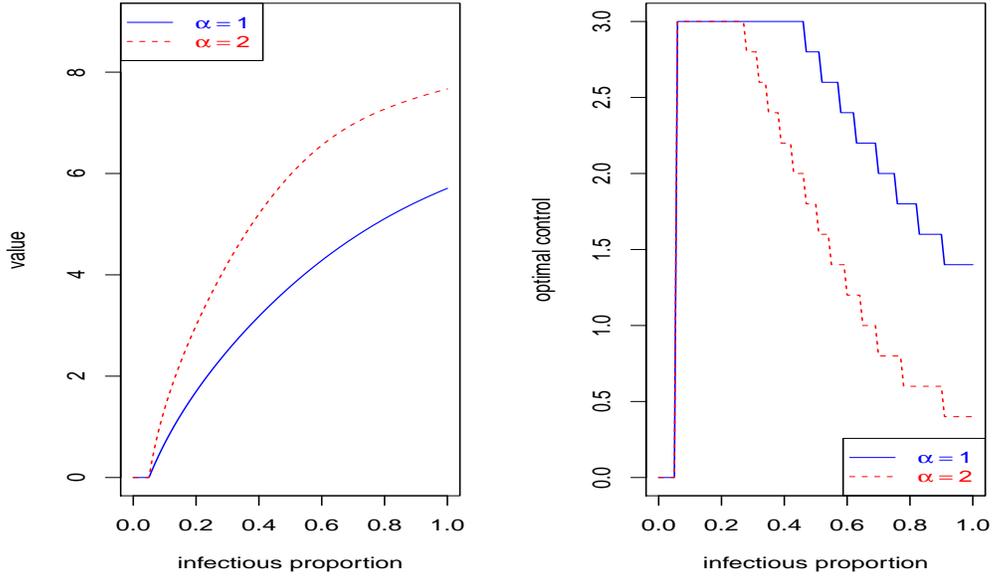}
		\caption{Value function (left) and optimal control (right) when $F(i, \ell, c)=1+\ell i+\ell ic^2$} \label{fig1A}
		\end{center}
\end{figure}

\begin{figure}[h!tb]
	\begin{center}
		\includegraphics[height=3.5in,width=5.5in]{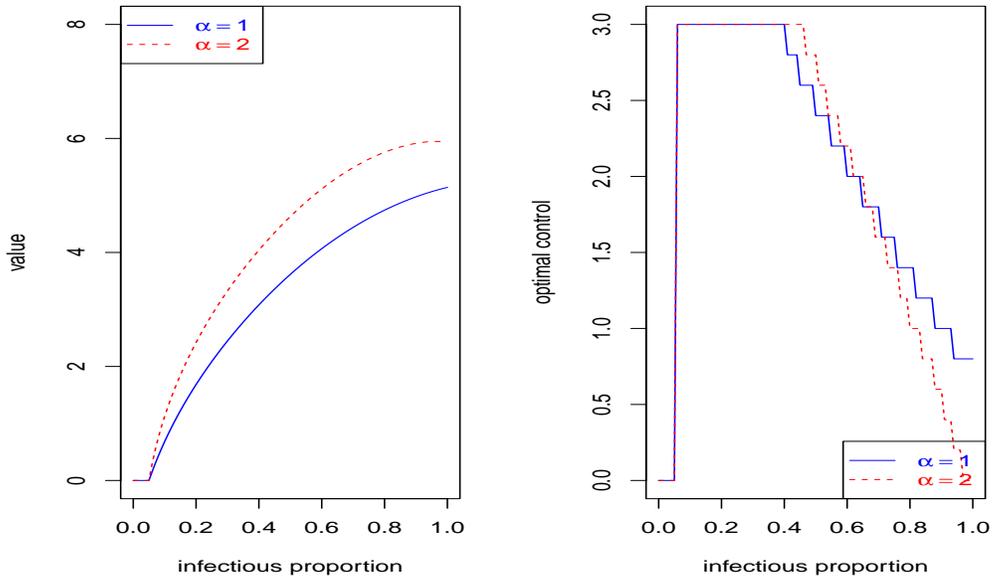}
		\caption{Value function (left) and optimal control (right) when $F(i, \ell, c)=1+2\lambda_\ell (1-i)i + ic^2$} \label{fig1B}
		\end{center}
\end{figure}

\begin{figure}[h!tb]
	\begin{center}
\includegraphics[height=3.5in,width=5.5in]{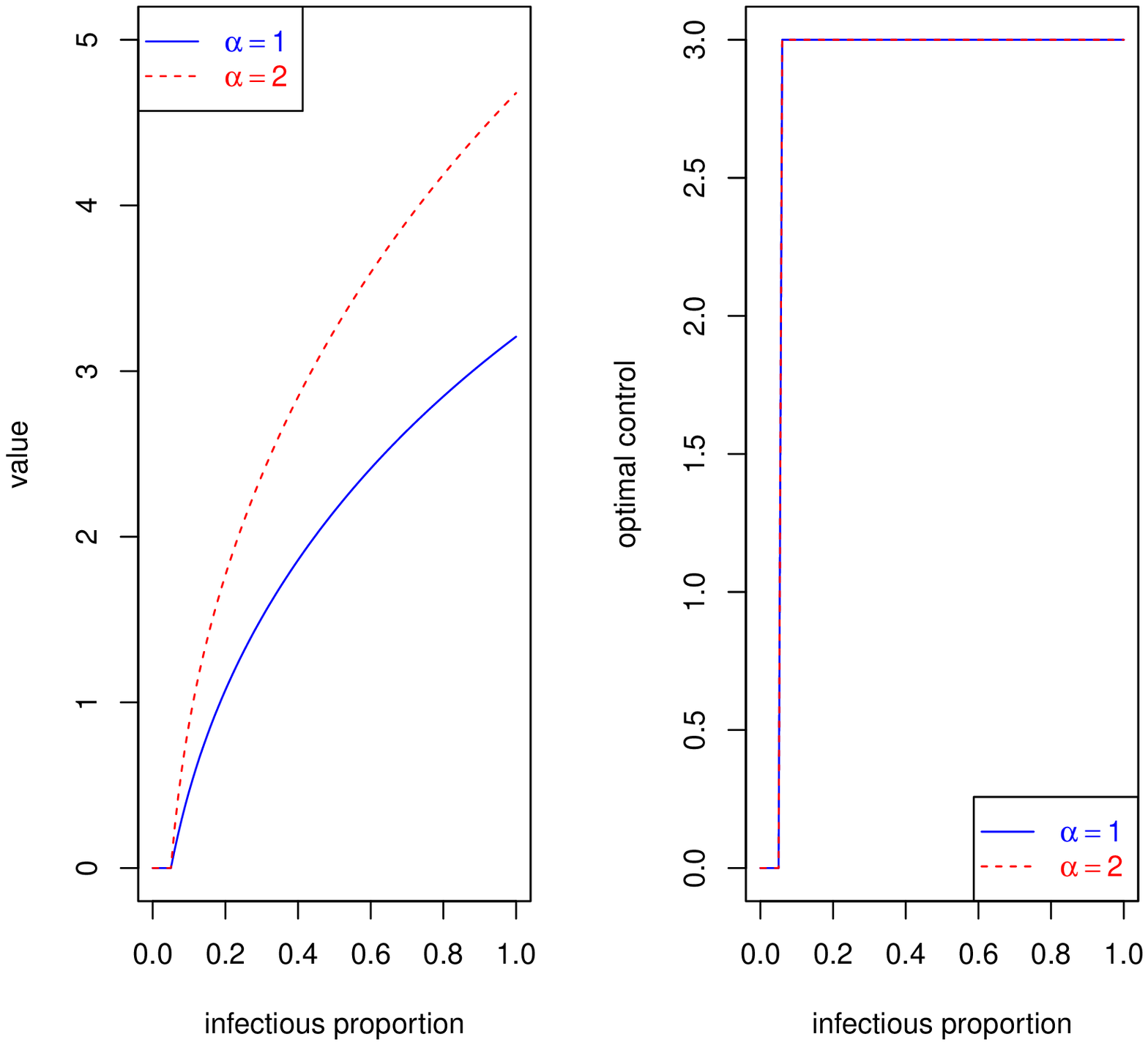}
		\caption{Value function (left) and optimal control (right) when $F(i, \ell, c)=1+\ell i+ \ell ic$} \label{fig1C}
		\end{center}
\end{figure}

In the second experiment, we consider the cost function
\beq{2e:ee2} F(i, \ell, c)=1+2\lambda_\ell (1-i)i + ic^2,\eeq
aiming at minimizing the sum of the eradication time and the total epidemic size. The value function and optimal control is shown in Figure \ref{fig1B}. The optimal control in each regime has a similar shape as that in the preceding experiment. Although $F(i, 2, c)>F(i, 1, c)$ for any $(i, c) \in (0, 1)\times \mathcal{U}$, the optimal control for $i\in (0.4, 0.6)$ in regime 2 is higher than that in regime 1.

In the third experiment, we consider the cost function
\beq{2e:ee} F(i, \ell, c)=1+\ell i+\ell ic.\eeq
The value function and optimal control is shown in Figure \ref{fig1C}.
The results in Figure \ref{fig1C} tells us that the control cost is small enough so that we should apply the maximum possible control in any regime.

}

\end{exm}

\begin{figure}[h!tb]
 	\centering		\subfloat{{\includegraphics[width=7.75cm]{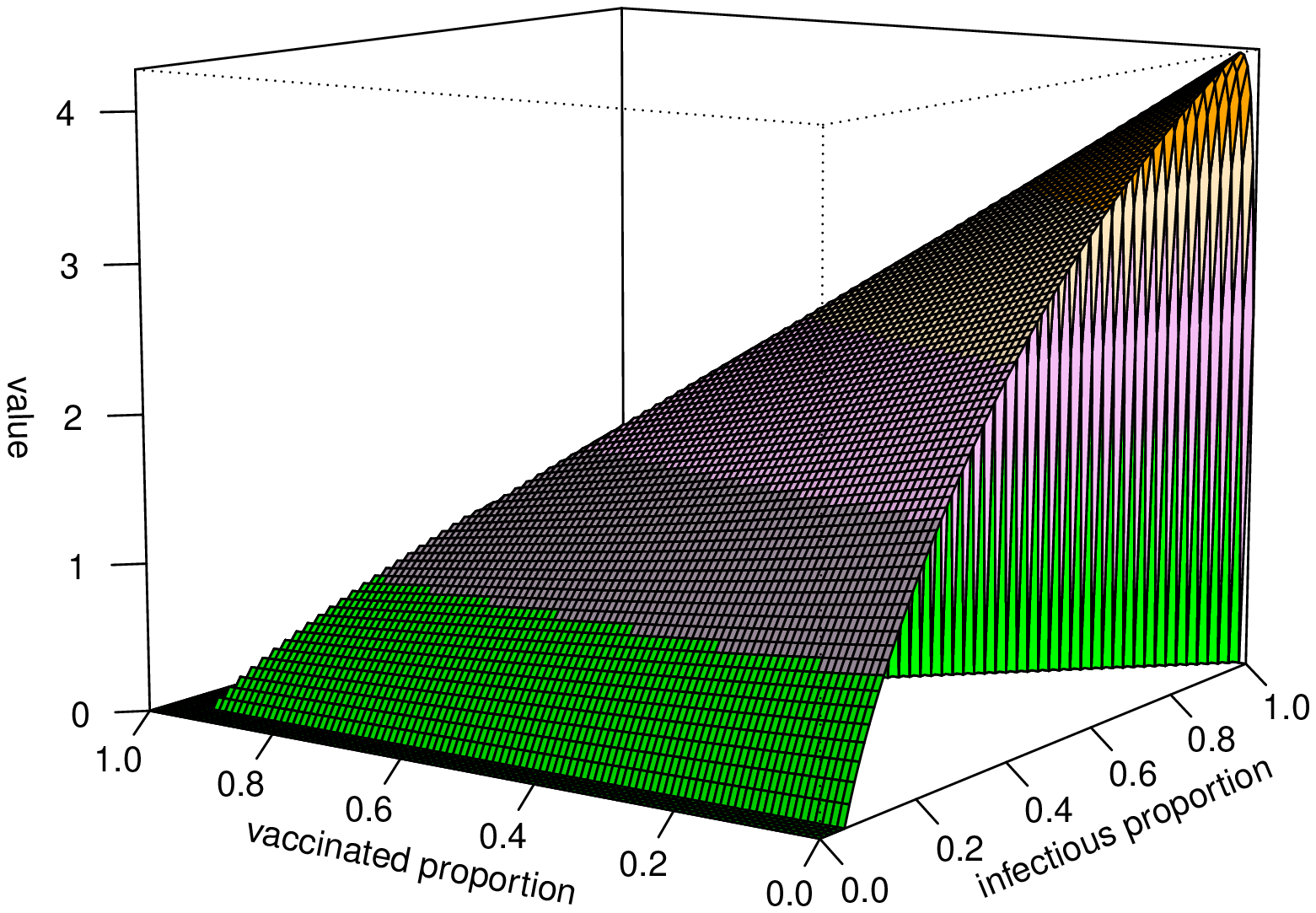} }}%
			\quad
\subfloat{{\includegraphics[width=7.75cm]{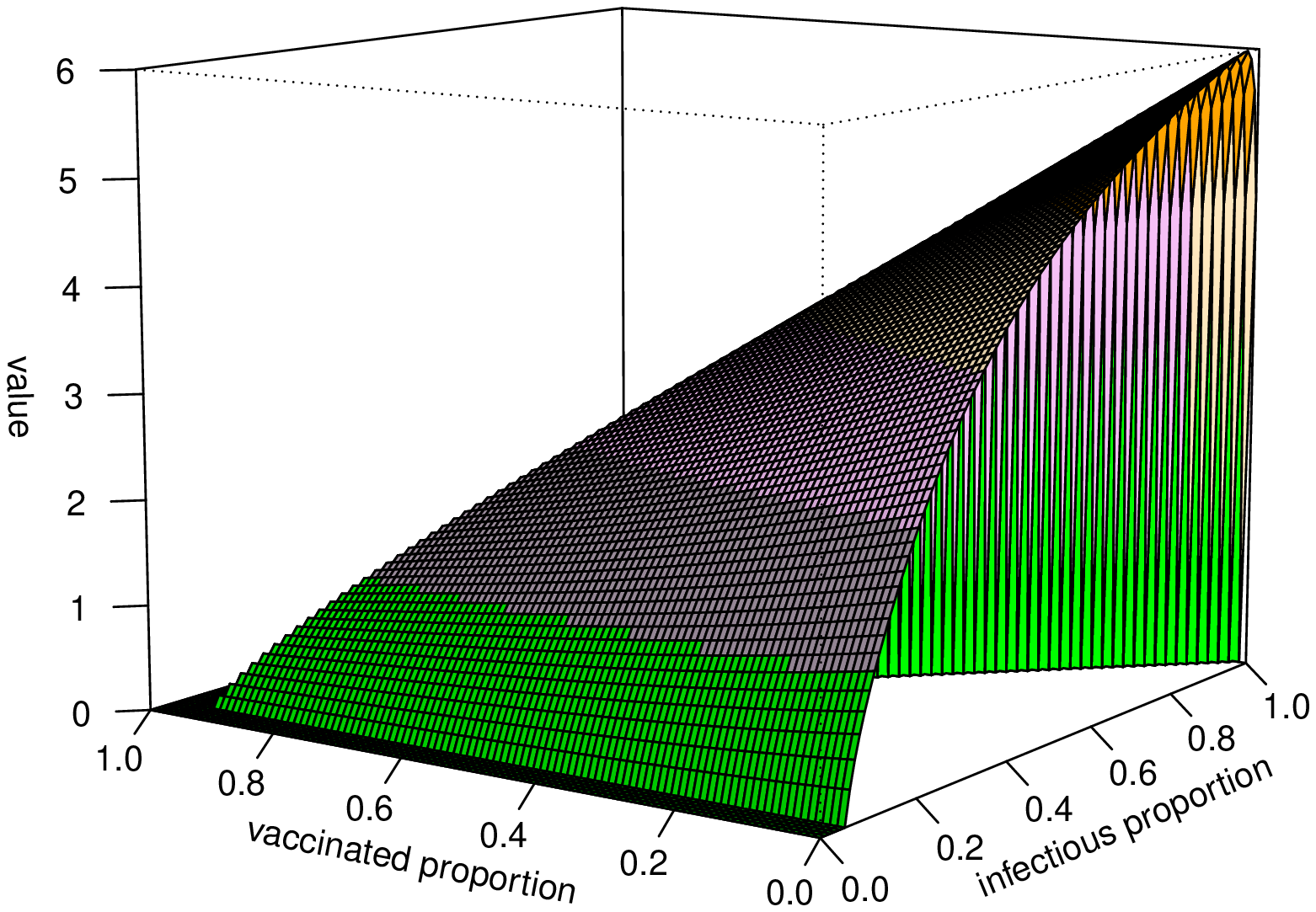} }}%
	\caption{The value function in regime 1 (left) and in regime 2 (right)}
	\label{fig2v}
\end{figure}

\begin{figure}[h!tb]
 	\centering		\subfloat{{\includegraphics[width=7.75cm]{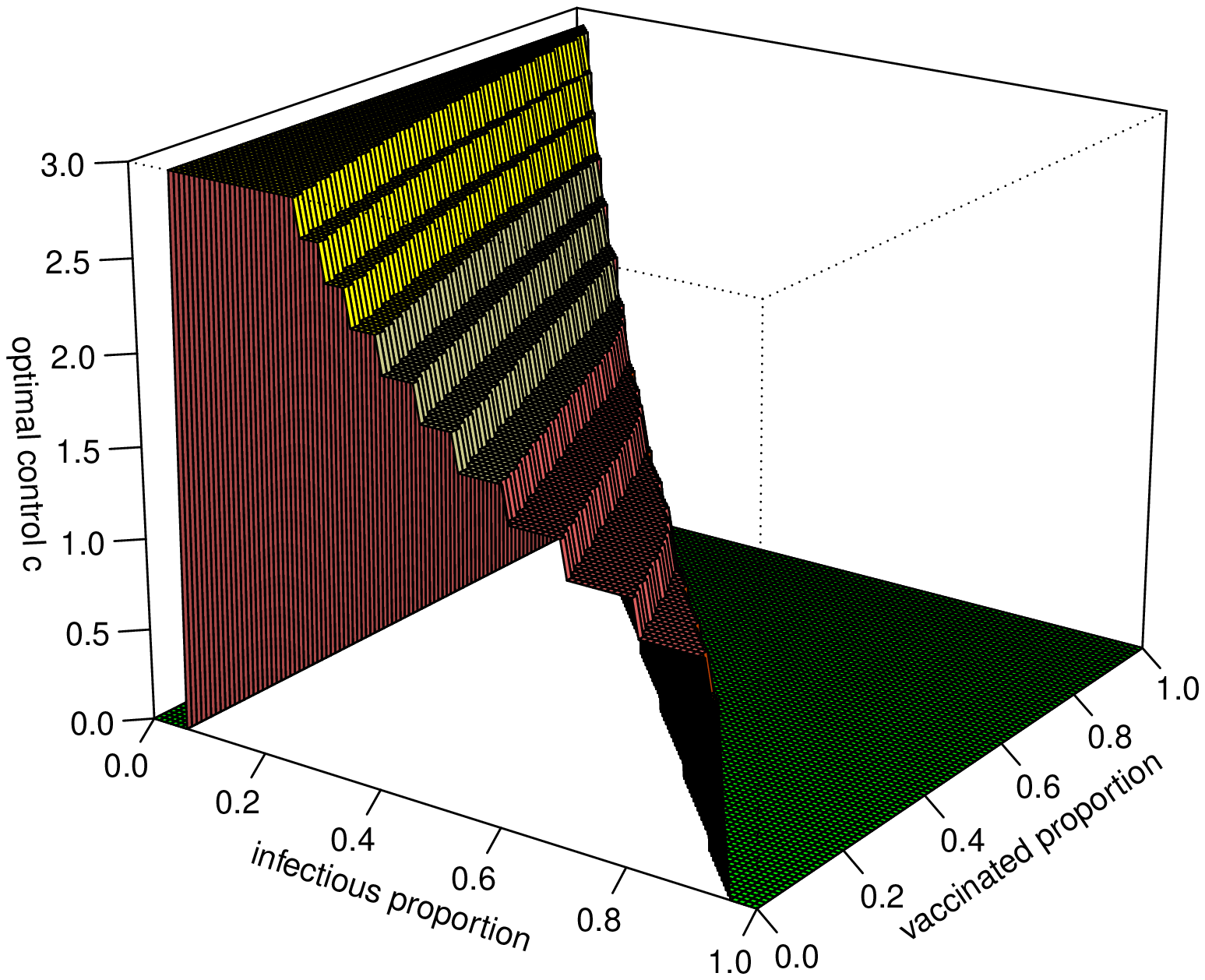} }}%
			\quad
\subfloat{{\includegraphics[width=7.75cm]{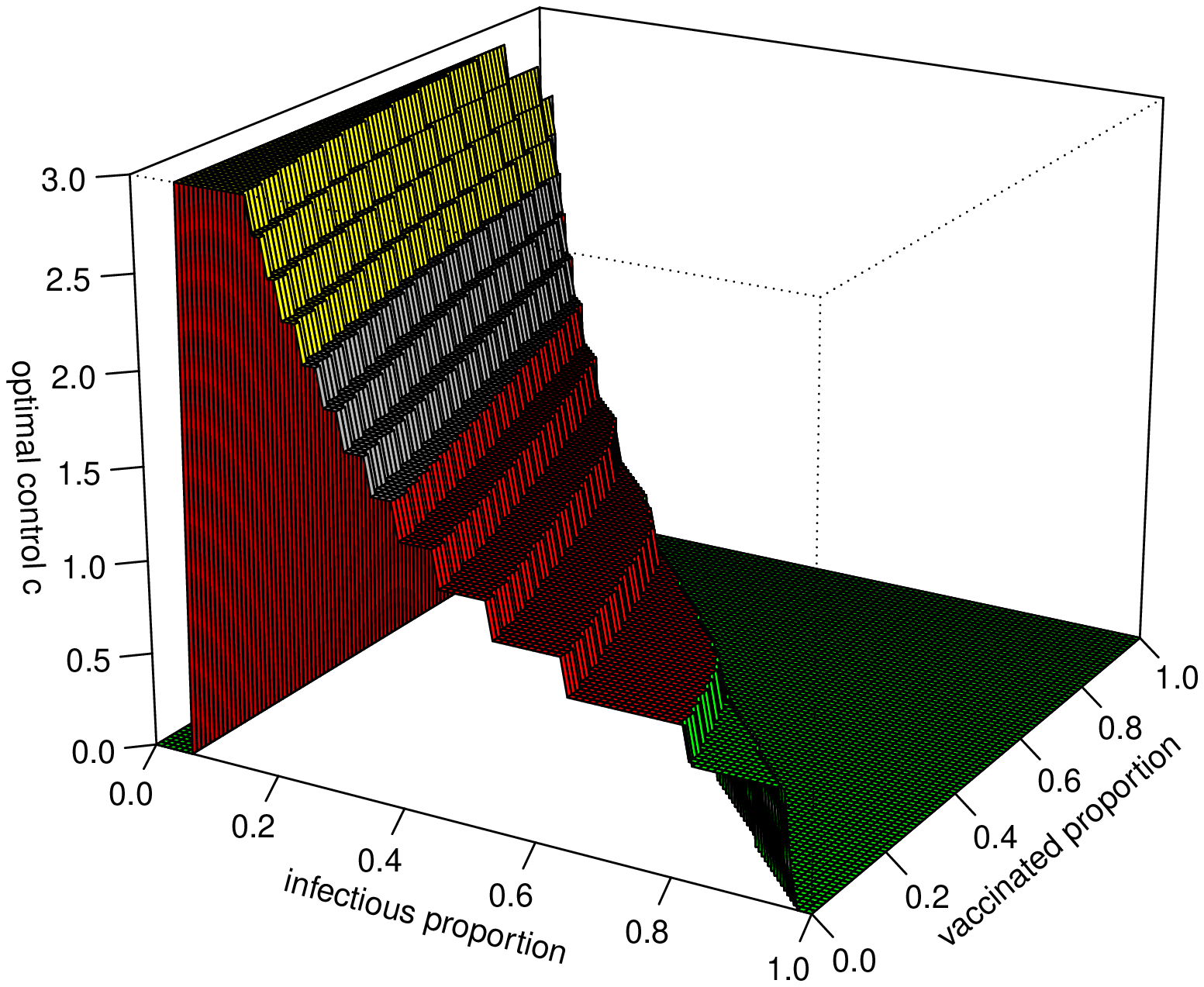} }}%
	\caption{The optimal control $C\cd$ in regime 1 (left) and in regime 2 (right)}
	\label{fig2c}
\end{figure}

\begin{exm} \label{ex3}
{\rm
We consider the model with vaccination given by \eqref{ee.3}; that is,
\beq{ee.5x}
\begin{cases}
{d I(t)} &= \big[ \lambda_{\al(t)} S(t)I(t) - \big(\mu_{\al(t)} +\gamma_{\al(t)} +C(t)\big)I(t) \big]dt + \sigma_{\al(t)} S(t)I(t)dw(t),\\
{d \V(t)} &=\big[ \mu_{\al(t)}q(t) +p(t)S(t) - (\mu_{\al(t)}+\e_{\al(t)}) \V(t) \big]dt.
\end{cases}
\eeq
We use the same parameters and the control set as in the preceding exam. In addition,
\begin{equation*}
\barray
\aad F(i, \ell, c, p, q) = 1 +\ell i +2 \ell ic^2 +(0.1)p + (0.1)q, \\
\aad \mathcal{V}_p = \mathcal{V}_q = \{0.2k: 0\le k\le 4\}.
\earray
\end{equation*}
We take the initial control
$$C_0(i, \vv,\ell)\equiv \max \mathcal{U}, \quad
p_0(i, \vv, \ell)\equiv \max \mathcal{V}_p, \quad
q_0(i, \vv, \ell)\equiv \max \mathcal{V}_q,$$
and set the initial value
 $\V^h_0(i, \vv, \ell)\equiv \int_0^\infty e^{-\delta t} F(1, 2, C_0, p_0, q_0)dt$ for $(i, \vv,\ell)\in T_h$.
We  outline how to find the sequence of costs of $V_n(\cdot)$ as follows.
For each $(i, \vv,\ell)\in T_h$, and control $(c, p, q)\in \mathcal{U}\times \mathcal{V}_p\times \mathcal{V}_q$,
 we compute
 \begin{equation*}
 \barray
V_{n+1}^h(i, \vv, \ell \,|\, c, p, q) \ad =e^{-\delta \Delta t^h(i, \vv, \ell, c, p, q) }\sum\limits_{(i', v, \ell')\in T_h} V^{h}_n (i', \vv',\ell'
) p^h \big((i, \vv, \ell), (i', \vv', \ell') | c, p, q \big)\\
\ad \qquad\qquad + F(i, \ell, c, p, q)\Delta t^h(i, \vv, \ell, c, p, q).
\earray
\end{equation*}
Then we choose the control  and record an improved value
$V^h_{n+1}(i, \vv, \ell)$ by
$$\big(C^h_{n+1}(i, \vv, \ell),  p^h_{n+1}(i, \vv, \ell), q^h_{n+1}(i, \vv, \ell)\big)= {\rm arcmin}_{(c, p, q)\in \mathcal{U}\times \mathcal{V}_p\times \mathcal{V}_q} V_{n+1}^h(i, \vv, \ell \,|\, c, p, q)$$
and
$$ V^h_{n+1} (i, \vv, \ell) = V^h_{n+1}\big(i, \vv, \ell\,|\, C^h_{n+1}(i, \vv, \ell),  p^h_{n+1}(i, \vv, \ell), q^h_{n+1}(i, \vv, \ell)\big).$$
The iterations stop as soon as the
increment
$V^h_{n+1}\cd-V^h_n\cd$
reaches the tolerance level. Intuitively, the shape of the optimal control will depend on how large the coefficients of the cost function are.
The value function and optimal control are displayed
in Figures  \ref{fig2v}, \ref{fig2c}, \ref{fig2p}, and \ref{fig2q}.

It can be seen that with the presence of an effective vaccination, the value function shown in Figure \ref{fig2v} is much lower than that of
Figure \ref{fig1A} in
Example 1. Moreover, Figure \ref{fig2c} reveals that one should use a lower control $C\cd$ compared to the case with no vaccination. Figures \ref{fig2p} and \ref{fig2q} provide the optimal vaccination plan. It appears that the optimal vaccination plan depends on the regime of the environment. The numerical studies emphasize the importance of effective vaccines with low costs in controlling epidemics.

\begin{figure}[h!tb]
 	\centering		\subfloat{{\includegraphics[width=7.75cm]{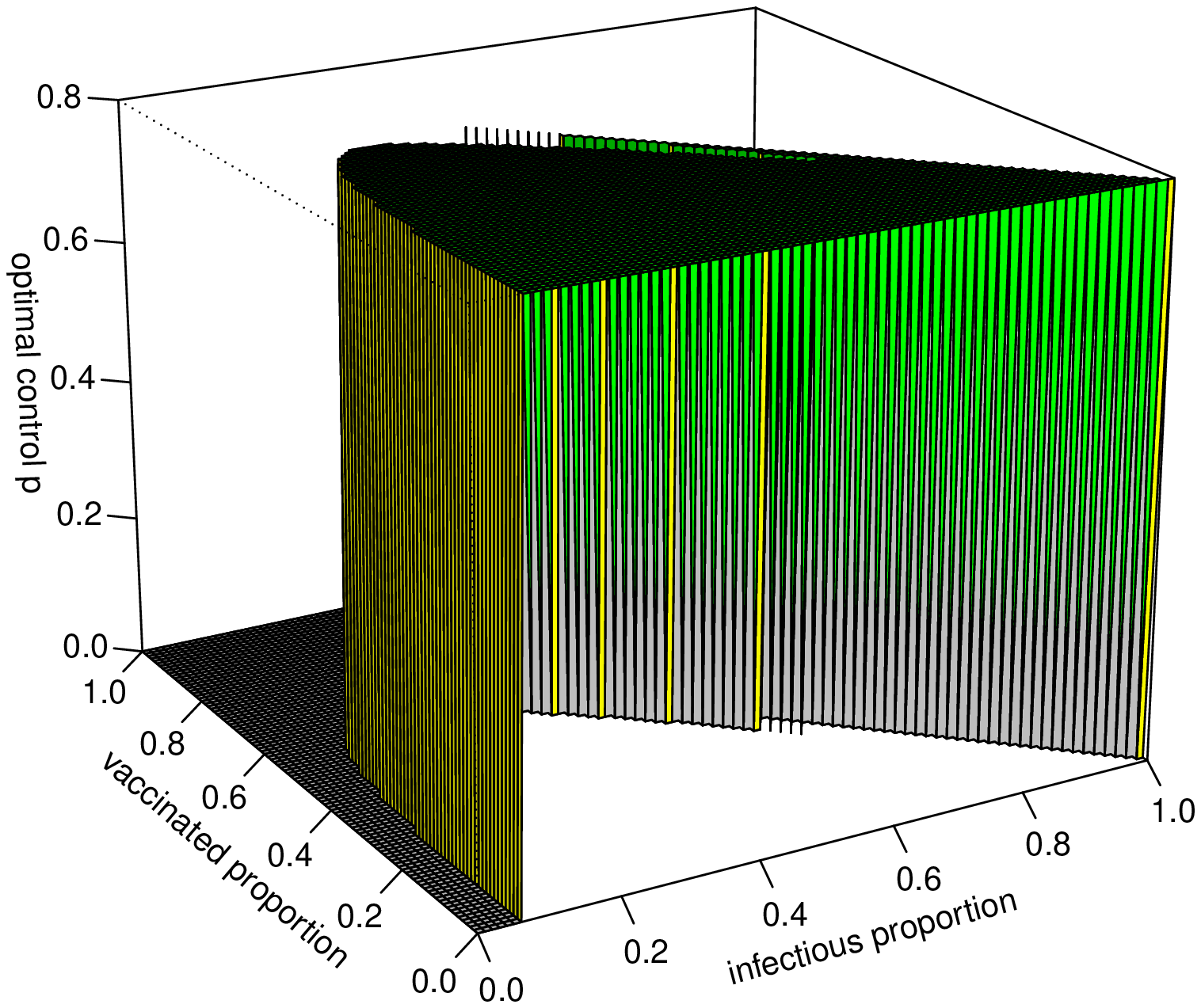} }}%
			\quad
\subfloat{{\includegraphics[width=7.75cm]{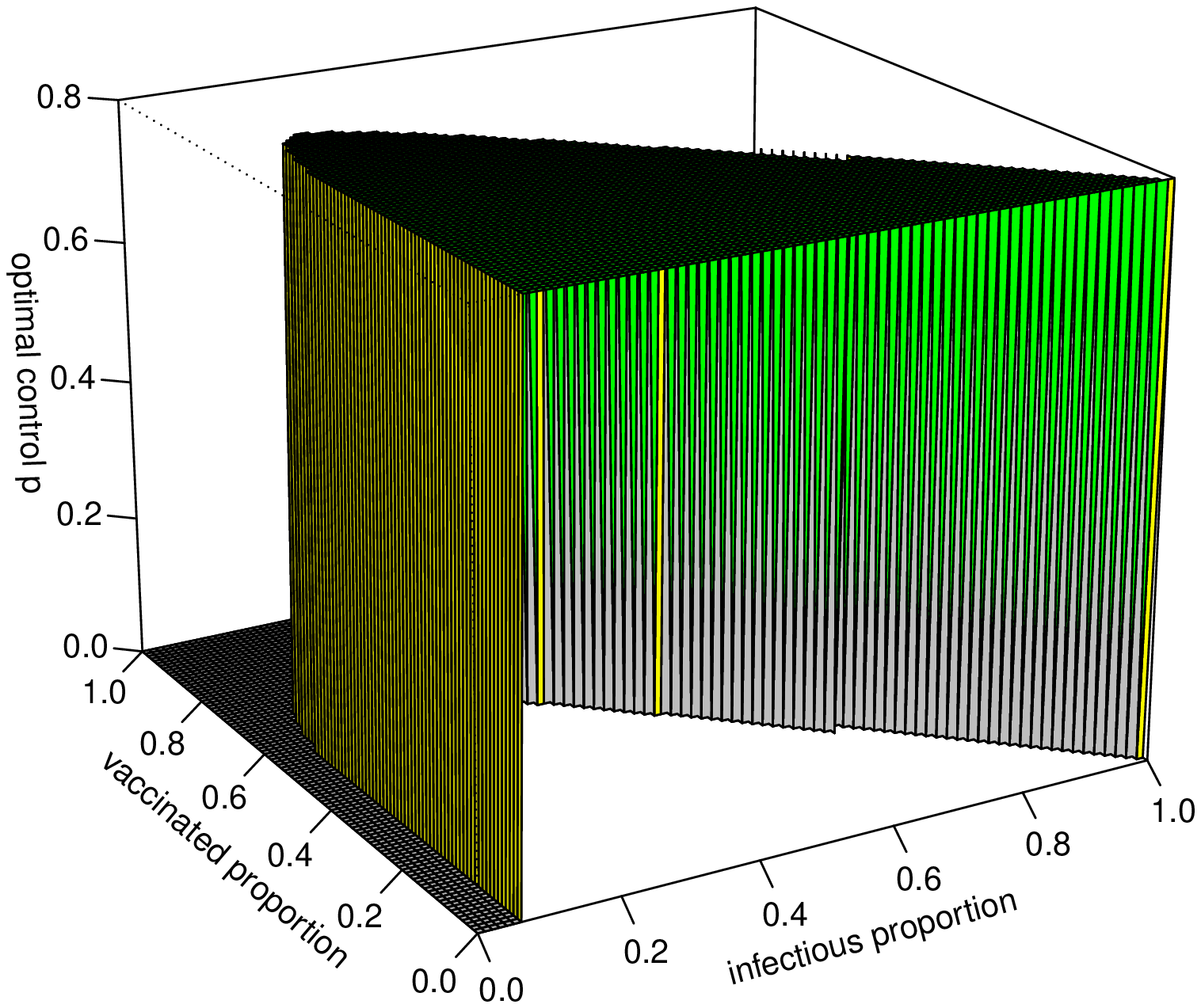} }}%
	\caption{The optimal control $p\cd$ in regime 1 (left) and in regime 2 (right)}
	\label{fig2p}
\end{figure}

\begin{figure}[h!tb]
 	\centering		\subfloat{{\includegraphics[width=7.75cm]{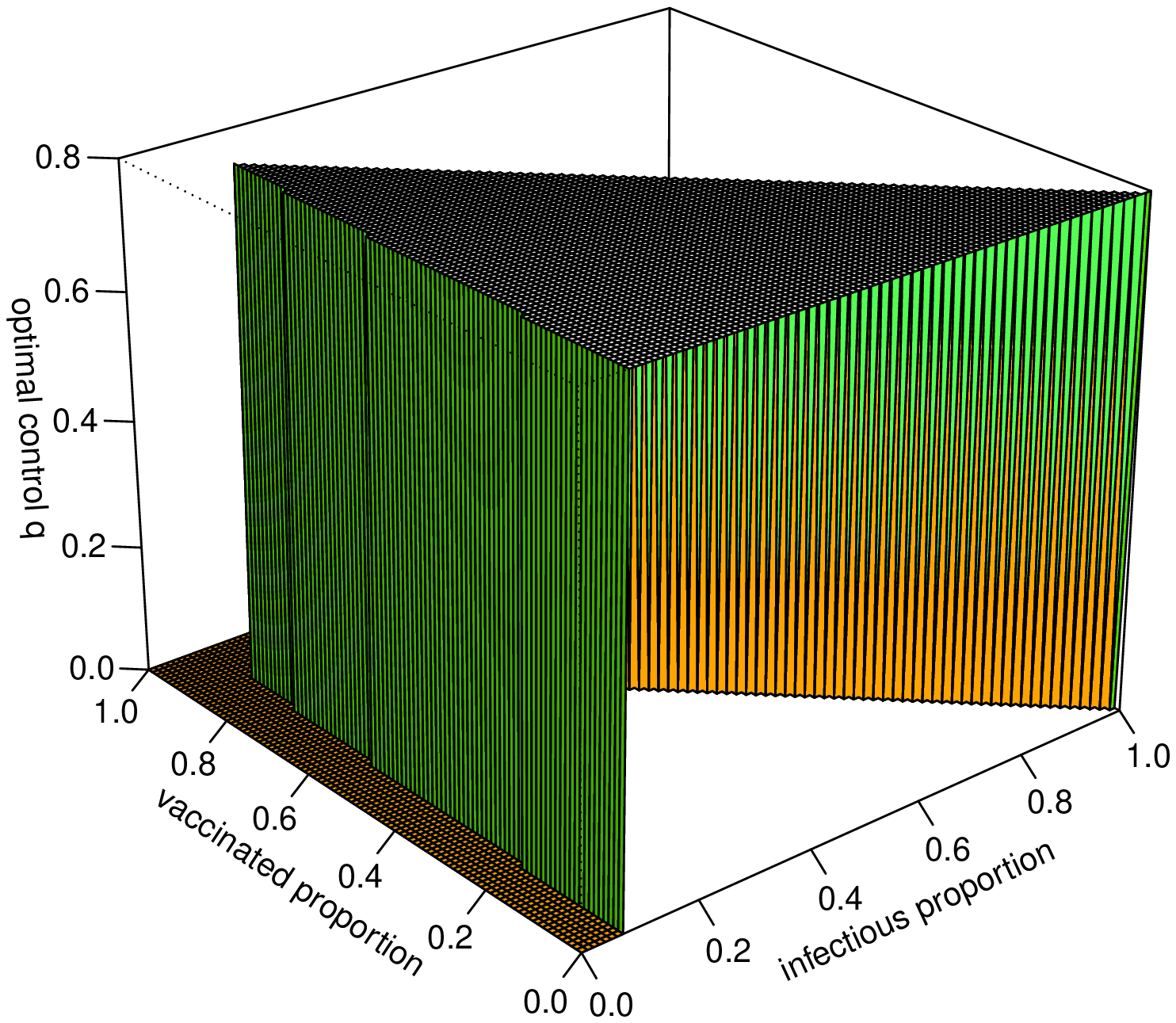} }}%
			\quad
\subfloat{{\includegraphics[width=7.75cm]{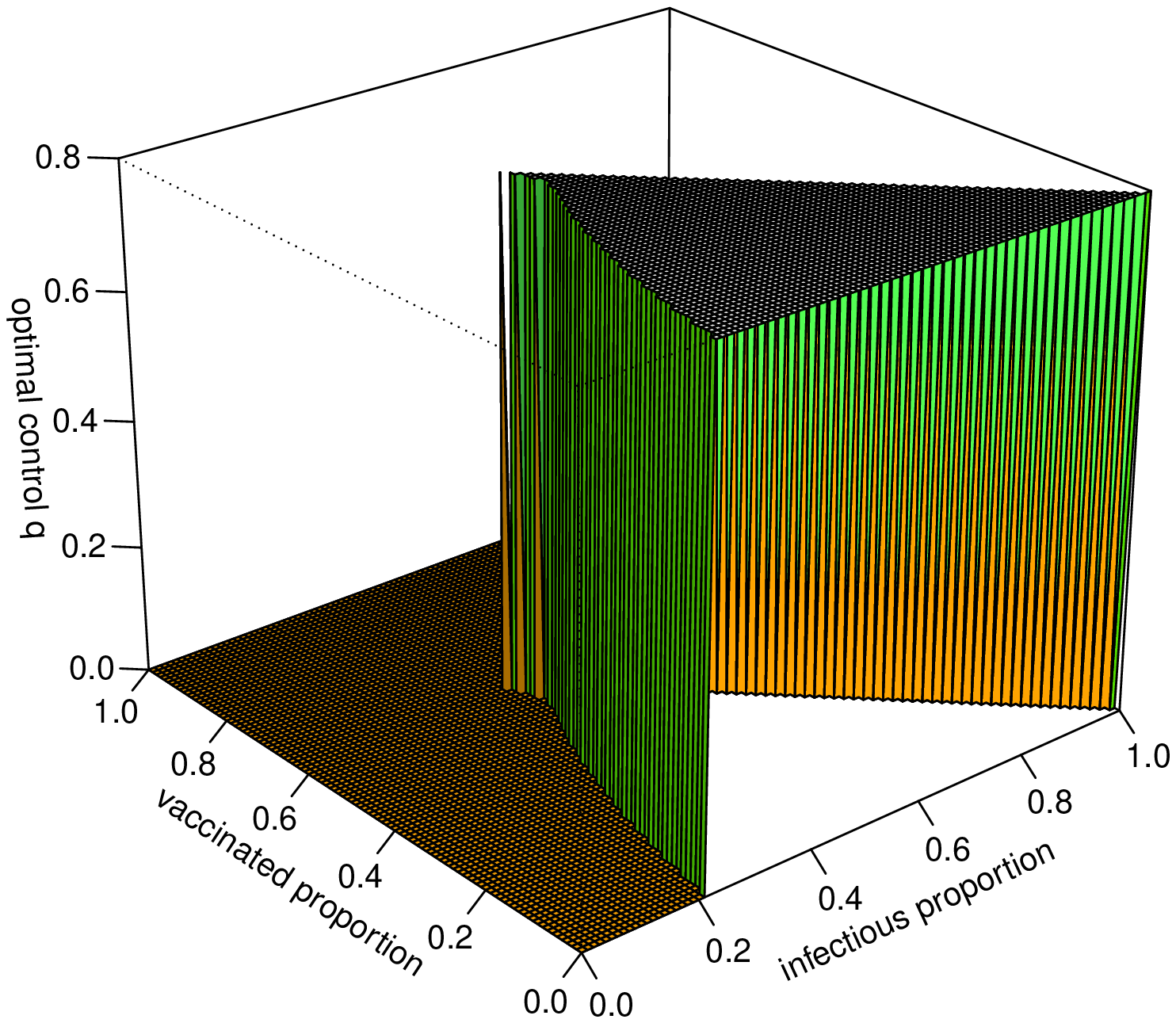} }}%
	\caption{The optimal control $q\cd$ in regime 1 (left) and in regime 2 (right)}
	\label{fig2q}
\end{figure}

}
\end{exm}

\section{Further Remarks}\label{sec:fur}
This paper focused on numerical methods for optimal control of  SIS epidemic models. We considered two models incorporating treatments, isolation, and vaccination. Using the Markov chain approximation method, we are able to treat a hybrid model affected by two types of environmental fluctuations with a general cost consideration.
The convergence of the algorithms was proved. Several numerical examples were  used to demonstrate the performance of our algorithm.
Some interesting questions deserve further investigation. One can study a set of mixed controls consisting of vaccination, isolation, and social distancing, etc. One can also apply the approach in this work to investigate a variety of control problems of epidemic models, which becomes an urgent issue due to
COVID-19.
Thus, the techniques and the simulation study might be of interests
to researchers in various disciplines.

\appendix

\section{Proofs of results}

\para{Proof of \thmref{thm:one}.}
The proof is a special case of that of Theorem \ref{thm:two} given below  when we take $\V(0)=p(t)=q(t)=0$ for any $t\ge 0$ with probability one. $\Box$

\

\para{Proof of \thmref{thm:two}.}
Since the coefficients of \eqref{ee.3} are locally Lipschitz continuous,
there is a unique local solution $(I(t), \V(t), \al(t))$ on $t\in [0,
\zeta)$, where $\zeta$ is the explosion time.
Now, by \eqref{ee.3}, $S+I+\V$ satisfies
 \begin{equation*}
\begin{cases}
& d (S+I+\V)(t) = \mu_{\al(t)} \big( 1- S(t) - I(t) - \V(t) \big)dt,\\
& S(0) + I(0) + \V(0) = 1,
\end{cases}
\end{equation*}
which implies
\beq{ee.5} S(t) + I(t) + \V(t) =1 \quad \text{for} \quad t\in [0, \zeta).
\eeq
Let $k_0$ be a sufficiently large positive integer
such that $s_0, i_0
\in\({1}/{k_0}, 1\)$. For each $k\ge k_0$, we define \beq{c1.2}
\tau_k:=\inf\left\{t\ge 0: I(t)\le \frac{1}{k} \quad \text{or}\quad S(t) \le \frac{1}{k} \right\}.\eeq Clearly the sequence $\{\tau_k\}$
is monotonically increasing. Let $\tau_\infty:=\lim\limits_{l\to
\infty}\tau_k$. Then $\tau_\infty\le \zeta$. It suffices to show
that $\tau_\infty=\infty$ with probability one. If this were false, there would
exist a $T>0$ and $\e>0$ such that $P\{\tau_\infty\le T\}> \e$.
Therefore we can find some $k_1\ge k_0$ such that \beq{c2}
P\{\tau_k\le T\}>\e \quad \text{for}\quad k\ge k_1.\eeq
We fix $\omega\in \Omega$. For $t\in [0, \tau_\infty)$, by \eqref{ee.5} and \eqref{c1.2}, $I(t)\in (0, 1]$.
We have from the third equation in \eqref{ee.3} that
\beq{ex}
\barray
\ad \V(t)\le \vv_0+ \int_0^t \Big( \mu_{\al(u)}q_0 + p_0 - \big[ \mu_{\al(u)} +\e_{\al(u)} + p_0\big]\V(u)\Big)du\le \vv_0+K_1 \int_0^t \big(1-\V(u)\big)du,\\
\ad \V(t)\ge \vv_0+\int_0^t \Big( \mu_{\al(u)}q_0- \big[ \mu_{\al(u)} +\e_{\al(u)} + p_0\big]\V(u)\Big)du \ge \vv_0 - K_2\int_0^t \V(u)du,
\earray
\eeq
where $K_1$ and $K_2$ are nonnegative numbers.
Since $\vv_0\in [0, 1)$ it follows from \eqref{ex} that $\V(t)\in [0, 1)$ for any $t\in [0, \tau_\infty)$. Combining this fact and \eqref{ee.5}, we have
$$(S(t), I(t), \V(t))\in (0, 1)^2\times [0, 1) \quad \text{for} \quad t\in [0, \tau_\infty) \quad \text{with probability one}.$$
To proceed, we
consider the function
$$\Phi(s, i)=\frac{1}{i} +\dfrac{1}{s}, \quad (s, i)\in (0, 1)^2,$$
and define
\begin{equation*}
\barray
\L \Phi(s, i, \vv,\ell) \ad= \big( \mu_\ell(1-q_0) -\lambda_\ell si - (\mu_\ell + p_0)s + \gamma_\ell i +\e_\ell \vv \big) \dfrac{\partial \Phi}{\partial s} (s, i) \\
\ad \qquad +\big( \lambda_\ell si - ( \mu_\ell +\ga_\ell+c_0)i  \big)\dfrac{\partial \Phi}{\partial i} (s, i)+ \frac{1}{2} \dfrac{\partial^2 \Phi}{\partial s^2} \sg^2_\ell s^2i^2+\frac{1}{2} \dfrac{\partial^2 \Phi}{\partial i^2}\sg^2_\ell s^2i^2.
\earray
\end{equation*}
Note that for our study, it suffices to consider $\Phi$ to be independent of the switching states.
  We have $\Phi(s, i)> 0$
 and
\begin{equation*}
\barray
\ad \dfrac{\partial \Phi}{\partial s} (s, i)=-\dfrac{1}{s^2} , \quad \dfrac{\partial \Phi}{\partial i} (s, i)=-\dfrac{1}{i^2},\\
\ad \dfrac{\partial^2 \Phi}{\partial s^2} (s, i)= \dfrac{2}{s^3}, \quad \dfrac{\partial^2 \Phi}{\partial i^2} (s, i)= \dfrac{2}{i^3}
\earray
\end{equation*}
It follows that \beq{c3} \barray  \L \Phi(s, i, \vv,\ell) \ad= -\big( \mu_\ell(1-q_0) -\lambda_\ell si - (\mu_\ell + p_0)s + \gamma_\ell i +\e_\ell \vv \big) \dfrac{1}{s^2} \\
\ad \qquad -\big( \lambda_\ell si - ( \mu_\ell +\ga_\ell+c_0)i  \big)\dfrac{1}{ i^2} + \frac{1}{2} \dfrac{2}{ s^3} (s, i)\sg^2_\ell s^2i^2+\frac{1}{2} \dfrac{2}{i^3} (s, i)\sg^2_\ell s^2i^2.
\\
\ad \le K\Phi(s, i),\earray \eeq for some positive constant $K$ independent of $(s, i, \vv,\ell)$.
By Dynkin's formla,
we obtain for $k\ge k_1$ and $t\in [0, T]$ that
\begin{equation*}
\barray
\E \Phi\(S(\tau_k\wedge t), I(\tau_k\wedge t)
\)\ad =\Phi\(s_0, i_0\)+\E\int\limits_0^{\tau_k\wedge t} \L
\Phi\(S(u), I(u),\V(u),\al(u)\)du\\
\ad \le \Phi\(s_0,i_0\)+K \E\int\limits_0^{t}
\Phi\(S(u\wedge \tau_k), I(u\wedge \tau_k)\)du.
\earray
\end{equation*}
 The Gronwall inequality yields that
\beq{ha}\E \Phi\(S(s\wedge \tau_k), I(\tau_k\wedge T)\)\le \Phi(s_0, i_0)e^{KT}.\eeq
 Note that for each $\omega\in \{\tau_k\le
T\}$,
$I(\tau_k(\omega))\le \frac{1}{k}$ or
$S(\tau_k(\omega))\le \frac{1}{k}$.
 It follows from the definition of $\Phi(\cdot, \cdot)$ that
$\Phi\big( S(\tau_k(\omega)), I(\tau_k(\omega))\big)\ge k$. In view of
\eqref{c2} and \eqref{ha}, we obtain
$$\Phi(s_0, i_0)e^{KT}\ge k\e,$$ leading to a
contradiction as $k\to \infty$. Thus, $\zeta = \tau_\infty=\infty$ with probability one.
The conclusion follows. $\Box$

The proofs below are motivated by \cite{Kushner92, Song2006}.
Let $D[0, \infty)$ denote the space of functions that are right continuous and have left-hand limits endowed with the Skorohod topology. All the weak convergence analysis
will be on this space or its $k$-fold products $D^k[0, \infty)$ for appropriate $k$. We will
provide a sketch of the proofs and refer to the well-known references for the details.

\

\para{Proof of \thmref{thm:thm}.}
(a)	The tightness of $\{\al^h\cd\}$ is obvious by Lemma \ref{lem:1}. For other components, we use the tightness criteria in \cite[p. 47]{Kushner84}. Specifically,  a sufficient condition for tightness of a sequence of processes  $\zeta^h\cd$ with paths in $D^k[0, \infty)$ is that for any $T_0, \rho
		\in (0, \infty)$,
		\bea
		\ad  \E_t^h\big|\zeta^h(t+s)-\zeta^h(t)\big|^2\le \E^h_t \gamma(h, \rho) \quad \text{for all}\quad s\in [0, \rho], \quad t\le T_0,\\
		\ad \lim\limits_{\rho\to 0}\limsup\limits_{h\to 0} \E  \gamma(h, \rho) =0.
		\eea
		The detailed proof of the tightness is standard; see \cite{Song2006}.
		Thus, $H^h\cd$ is tight.
	As a result, a subsequence of $H^h\cd$ converges weakly to the limit $H\cd=\big(I\cd, \al\cd, B\cd, M\cd, m\cd, \tau\big)$. Since the sizes of jumps of $I^h\cd$, $B^h\cd$, and $M^h\cd$ go to zero as $h\to 0$, then $I\cd$, $B\cd$, and $M\cd$ have continuous paths with probability one.
	
  (b) follows from (a) and \eqref{e.4.15}.

   (c)
	Let ${\E}_t^h$ denote the expectation conditioned on ${\mathcal{F}}^h(t)$. By the definition of $M^{h}\cd$, it is an $\{{\mathcal{F}}^h(t)\}$-martingale with quadratic variation process
	 $$\int_0^t  a\big( I^h(u),  \al^h(u)\big) du +  \e^h_3(t),$$
	 where $\{ \e_3^h\cd\}$
  is an $ \{{\mathcal{F}}^h(t)\}$-adapted process satisfying $\lim\limits_{h\to 0} \sup\limits_{t\in [0, T_0]}\E| \e_3^h(t)|=0$  for $T_0\in (0, \infty).$ By the Burkholder-Gundy inequality, there is a positive constant $K$  such that
  $$\E | M^{h}(t)|^2\le K\E \Big| \int_0^t  a\big( I^h(u), \al^h(u)\big) d u +  \e^h_3(t)\Big|\le K(t+1),$$
  for any $h>0$. Thus, the family $\{  M^{h}(t): h>0 \}$ is uniformly integrable.
   For any $\rho>0$,
\beq{e:250}
\barray
\aad \E_t^h \big(M^{h}(t+\rho)-M^{h}(t)\big)=0,\\
\aad \E_t^h \big[{M}^{h}(t+\rho) - M^{h}(t)\big]\big[{M}^{h}(t+\rho) - M^{h}(t)\big]'
= \E_t^h \int_t^{t+\rho}  a(I^h(u), \al^h(u))du   +  \e^h_4(\rho),
\earray
\eeq
 where $\E|{\e}^h_4(\rho)|\to 0$ as $h\to 0$.
 To characterize $M^{ h}\cd$, let $r$ be an arbitrary integer, $t>0$, $\rho>0$ and $\{t_k: k\le r\}$ be such that $t_k\le t<t+\rho$ for each $k$.
Let $\Psi\cd$ be a real-valued and continuous function of its arguments with compact support. Then in view of \eqref{e:250}, we have
\beq{e:251}
\E\Psi({H}^h(t_k), k\le r)\big[ M^{h}(t+\rho)-M^{h}(t)\big]=0,
\eeq
and
\beq{e:2513}
\barray
\E\Psi({H}^h(t_k), k\le r)\Big(  \big[{M}^{h}(t+\rho) \ad - M^{h}(t)\big]\big[{M}^{h}(t+\rho) - M^{h}(t)\big]'\\
\ad \qquad
-\int_t^{t+\rho}  a(I^h(u), \al^h(u))du \Big)=\e^h_5(\rho),
\earray
\eeq
where $\E|{\e}_5^h(\rho)|\to 0$ as $h\to 0$.
By using the Skorohod representation, letting $h\to 0$ in \eqref{e:251}, we obtain
\beq{e:252}
\E\Psi({H}(t_k), k\le r)\big[ M(t+\rho)-M(t)\big]=0.
\eeq
Since $M\cd$  has continuous paths with probability one,
\eqref{e:252} implies that $M\cd$ is a continuous ${\mathcal{F}}\cd$-martingale. Moreover, \eqref{e:2513} gives us that
\beq{}
\barray
\E\Psi({H}^h(t_k), k\le r)\Big(  \big[{M}(t+\rho) \ad - M(t)\big]\big[{M}(t+\rho) - M(t)\big]'\\
\ad \qquad
-\int_t^{t+\rho}  a(I(u), \al(u))du \Big)=0,
\earray
\eeq
which gives us the quadratic variation of $M\cd$. The conclusion follows from the martingale representation theorem (see \cite[Theorem 3.4.2]{Ka88}).

   (d) follows immediately from the results in  (a), (b), and (c). $\Box$

\

\para{Proof of \thmref{thm:4.6} and \thmref{thm:4.6.3}.}
The proofs are
modifications of that for \cite[Theorem 7]{Song2006}. Hence, we omit the details for brevity. $\Box$

\end{document}